\documentclass[11pt, twoside, a4paper]{article}
\usepackage[utf8]{inputenc}
\usepackage{amsmath, amssymb, amsthm}            
\usepackage{amstext, amsfonts, a4}
\usepackage[english]{babel}
\usepackage{bbm}
\usepackage{graphicx}
\usepackage[margin=2.5cm]{geometry}
\usepackage[hidelinks]{hyperref}

\theoremstyle{plain}
	\newtheorem{theorem}{Theorem}[section]
	\newtheorem{proposition}[theorem]{Proposition}    
	\newtheorem{lemma}[theorem]{Lemma}          
	\newtheorem{corollary}[theorem]{Corollary}
	\newtheorem{conjecture}[theorem]{Conjecture}
	
\theoremstyle{definition}
	\newtheorem{definition}[theorem]{Definition}

\DeclareMathOperator{\di}{div}
\DeclareMathOperator{\diam}{diam}

 \def\NN{{\mathbb N}}  
  \def\RR{{\mathbb R}}

 \def\cB{{\mathcal B}} \def\cC{{\mathcal C}} 
\def\cD{{\mathcal D}}   
 \def\cH{{\mathcal H}}  
  \def\cL{{\mathcal L}} 
   
\def\cP{{\mathcal P}} \def\cQ{{\mathcal Q}}  
\def\cS{{\mathcal S}} \def\cT{{\mathcal T}}

  \def\mC{{\mathfrak C}} 
\def\mD{{\mathfrak D}}

\def\mS{{\mathfrak S}}

\begin{document}
\title{Vortex collapses for the Euler and Quasi-Geostrophic Models}
\author{\renewcommand{\thefootnote}{\arabic{footnote}}
Ludovic Godard-Cadillac\footnotemark[1]}
\footnotetext[1]{Universit\`a degli studi di Torino, Dipartimento di matematica, Via Carlo Alberto 10, Torino (Italia)\\
\emph{and} Nantes Université, Laboratoire de Mathématiques Jean Leray, 2 Chem. de la Houssinière, Nantes (France).\\
\textbf{keywords:} Fluid mechanics, Euler equations, Quasi-geostrophic equations, Point-vortices, Differential equations}
\date{\today}
\maketitle

\begin{abstract}
This article studies point-vortex models for the Euler and surface quasi-geostrophic equations. In the case of an inviscid fluid with planar motion, the point-vortex model gives account of dynamics where the vorticity profile is sharply concentrated around some points and approximated by Dirac masses. This article contains two main theorems and also smaller propositions with several links between each other. The first main result focuses on the Euler point-vortex model, and under the non-neutral cluster hypothesis we prove a convergence result. The second result is devoted to the generalization of a classical result by Marchioro and Pulvirenti concerning the improbability of collapses and the extension of this result to the quasi-geostrophic case.
\end{abstract}

\section{Introduction}

\subsection{The inviscid surface quasi-geostrophic equation}
We are interested in a model of point-vortices for the inviscid surface quasi-geostrophic equation
\begin{equation}\tag{SQG}\label{eq:SQG}
\left\{\begin{split}
&\partial_t\,\omega+v\cdot\nabla\omega=0,\\
&v=\nabla^\perp(-\Delta)^{-s}\omega,
\end{split}\right.
\end{equation}
where $v:\RR_+\times\RR^2\to\RR^2$ is the fluid velocity and $\omega:\RR_+\times\RR^2\to\RR$ is called the active scalar. The notation $\perp$ refers to the counterclockwise rotation of angle $\frac{\pi}{2}$. 
Typically, the surface quasi-geostrophic equation models the dynamic in a rotating frame of the potential temperature for a stratified fluid subject to Brunt-Väisälä oscillations. 
This is a standard model for geophysical fluids and it is intensively used for weather forecast and climatology.
For more details about this physical model, see e.g.~\cite{Pedlowsky_1987} or~\cite{Vallis_2006}.
Mathematically, the quasi-geostrophic equation has many properties in common  with the two-dimensional Euler equation written in terms of vorticity
\begin{equation}\tag{Euler 2D}\label{eq:Euler}
\left\{\begin{split}
&\partial_t\,\omega+v\cdot\nabla\omega=0,\\
&v=\nabla^\perp(-\Delta)^{-1}\omega.
\end{split}\right.
\end{equation}
The two-dimensional Euler equation can be seen as a particular case of the quasi-geostrophic equation where $s$ is equal to $1$. 
Local well-posedness of classical solutions for (SQG) was established in~\cite{Constantin_Majda_Tabak_1994}, where are also studied the analogies with the two and three-dimensional Euler equation. 
Solutions with arbitrary Sobolev growth were constructed in~\cite{Kiselev_Nazarov_2012} in a periodic setting. 
So far and contrarily to the two-dimensional Euler equations, establishing global well-posedness of classical solutions for (SQG) is an open problem. 
Note also that the global existence of weak solutions in $L^2(\RR^2)$ was established in~\cite{Resnick_1995}, but below a certain regularity threshold, these weak solutions show dissipative behaviors and non-uniqueness is possible~\cite{Buckmaster_Shkoller_Vicol_2016}.
Exhibiting global smooth solutions or patch solutions is a challenging issue as there is no equivalent of the Yudovitch theorem~\cite{Yudovich_1963}.
A first example was recently provided in~\cite{Castro_Cordoba_Gomez_2020} by developing a bifurcation argument from a specific radially symmetric function. 
The varia\-tional construction of an alternative example in the form 
of a smooth traveling-wave solution was completed 
in~\cite{Gravejat_Smets_2019, Godard-Cadillac_2019}. 
Corotating patch solutions with two patches~\cite{Hmidi_Mateu_2017} and $N$ patches forming an $N$-fold symmetrical pattern~\cite{Garcia_2020} were recently exhibited with bifurcation 
argument. The $\cC^1$ analogous of these solutions has also been investigated independently recently in~\cite{Ao_Davila_DelPinto_Musso_Wei_2020}.
Another recent independent result~\cite{Godard-Cadillac_Gravejat_Smets_2020} also build corotating solutions with $N$ patches using varia\-tional argument. In this last article, the desingularization of the associated point-vortex problem is achieved. 

The present work aims at developing the understanding of the links 
between the quasi-geostro\-phic equation and the two-dimensional 
Euler equation through the study of the point-vortex model.
In the case of the two-dimensional Euler equation, the 
point-vortex model 
is a system of differential equations for points on $\RR^2$ that 
approximates situations where the vorticity $\omega$ is highly 
concentrated around several points. 
In such a situation, it is more convenient to see the vorticity as 
being a sum of Dirac masses evolving in time. 
This model is widely studied in fluid mechanics of the plane. An extensive 
presentation of the main results on this system can be found 
at~\cite[Chap.\,4]{Marchioro_Pulvirenti_1994}, 
completed by~\cite{Marchioro_Pulvirenti_1984}. 
The desingularization problem which consists in a rigorous derivation of the point-vortex model is a classical issue
for the two-dimensional Euler equation~\cite{Smets_VanSchaftingen_2010}, to our 
knowledge it is still open for (SQG) vortices although recent
results exist~\cite{Geldhauser_Romito_2020, Rosenzweig_2020}.

This article generalizes several existing 
results known for the Euler vortices or extends known results for Euler to the quasi-geostrophic case.
The first proposition of this article is the generalization of a uniform bound result, Theorem 2.1 in \cite[Chap.\,4]{Marchioro_Pulvirenti_1994}. 
We prove that under the \emph{non-neutral cluster hypothesis} 
(defined hereafter) the vortices stays bounded in finite time and this bound does not depend on the singularity of the kernel nor on the initial position of the 
vortices. 
We also provide a uniform relative bound for a slight relaxation of the \emph{non-neutral cluster hypothesis}. 
In the second part of this work we prove that under the non-neutral cluster hypothesis the trajectories of the vortices for the Euler model are convergent in finite time even in the case of collapses. The quasi-geostrophic case is left open.
The third part of this work is devoted to the question of the 
improbability of collapses. 
This consists in studying the Lebesgue measure of the initial conditions leading to a collapse, which is expected to be equal to $0$. 
This question has been successfully answered by Theorem 2.2 in \cite[Chap.\,4]{Marchioro_Pulvirenti_1994} under the \emph{non-neutral cluster hypothesis} for the Euler point-vortices.
The extension to the quasi-geostrophic case was achieved by~\cite{Geldhauser_Romito_2020} in the case $1/2<s<1$.
We generalize this result to (SQG) for all $s\in(0,1]$ and we weaken the \emph{non-neutral cluster hypothesis} since we allow the total sum of the intensities of the vortices to be equal to $0$.

\subsection{Presentation of the point-vortex model}~\label{sec:presentation}
The point-vortex model on the plane $\RR^2$ consists in assuming that at time $t=0$ the vorticity can write as a sum of Dirac masses,
\begin{equation}\label{eq:vorticity Dirac}
\omega(t=0,x)=\sum_{i=1}^Na_i\delta_{x_i}.
\end{equation}
The points $x_i$ are the respective position of the vortices $a_i\delta_{x_i}$ and the coefficients $a_i\neq0$ are their intensity. 
The first equation in~\eqref{eq:SQG} or~\eqref{eq:Euler} is a transport equation on the vorticity $\omega$. 
It is expected that the Dirac masses initially located at $x_i$ are left unchanged but transported by the flow. 
Formally, if we solve the evolution equations~\eqref{eq:SQG} or~\eqref{eq:Euler} with initial datum~\eqref{eq:vorticity Dirac}, we obtain that the initial speed writes
\begin{equation}
\label{eq:vorticity Dirac in speed}
v(t=0,x)=\sum_{i=1}^Na_i\nabla^\perp_xG\big(|x_i-x|\big),
\end{equation}
where $G$ is the profile of the Green function of the fractional Laplace operator $(-\Delta)^s$ in the plane $\RR^2$. Here the parameter $s$ is chosen to be in $(0,1]$. In this case, The profiles $G:\RR_+^\ast\to\RR$ are given by:
\begin{equation}\label{def:Green functions}
G_1(r):=\frac{1}{2\pi}\log\Big(\frac{1}{r}\Big)\qquad\mathrm{and}\qquad G_s(r):=\frac{\Gamma(1 - s)}{2^{2s} \pi \Gamma(s)}\,\frac{1}{r^{2(1 - s)}},
\end{equation}
with $\Gamma$ the classical Gamma function. 
Nevertheless, a problem arises from the singularity of the speed in~\eqref{eq:vorticity Dirac in speed}.
Since the vorticity is concentrated in one point then it is usually assumed that it does not interact with itself but only with the vortices that are at a positive distance. 
We derive that the differential equation describing the evolution of the position of the vortices is given by
\begin{equation}\label{eq:vortex equation}
\frac{d}{dt}x_i(t)=\sum_{\substack{j=1\\j\neq i}}^Na_j\nabla^\perp G\big(|x_i(t)-x_j(t)|\big),
\end{equation}
where $\nabla G(|x|)$ is a shortcut for the gradient of the function $x\mapsto G(|x|)$. This useful notation is coming from the reference work~\cite[Chap.\,4]{Marchioro_Pulvirenti_1994}. In the other cases and to avoid ambiguity, the gradient of a function will be denoted using the variable in subscript: $\nabla_{x_i}$.
We are going to generalize some of the results in~\cite[Chap.\,4]{Marchioro_Pulvirenti_1994} to more general kernel profiles $G$ that include the quasi-geostro\-phic case. 
In the sequel, the function $G:\RR_+^\ast\to\RR$ is assumed to be chosen such that~\eqref{eq:vortex equation} satisfy the hypothesis of the Cauchy-Lipschitz theorem (also known as Picard-Lindelöf theorem) as long as the distances between vortices remain positive.
\begin{definition}[Set of collapses]\label{defi:collisions}
The set of initial datum such that two or more vortices collapse on the interval of time $[0,T)$ is defined by
\begin{equation}\label{def:collisions T}
\mC_T:=\Big\{X\in\RR^{2N}\;:\;\exists\;T_X\in[0,T),\;\;\liminf\limits_{t\to T_X}\;\min\limits_{i\neq j}\;|x_i(t)-x_j(t)|=0\Big\}.
\end{equation}
We then set
\begin{equation}
\mC:=\bigcup_{T=1}^{+\infty}\mC_T.
\end{equation}
The set $\mC$ is called the \emph{set of collapses} and $T_X$ is the time of collapse associated to the initial datum $X\in\mC.$ Note that these sets depend on the choice of the kernel $G$.
\end{definition}
The point-vortex differential equation~\eqref{eq:vortex equation} is well-defined for all initial datum $X\in\RR^{2N}\setminus\mC$. 
If we restrict the analysis to a bounded interval of time $[0,T)$ then it is well-defined for all initial datum $X\in\RR^{2N}\setminus\mC_T$, which eventually allows us to study a possible vortex collapse at time $t=T.$

Concerning the point-vortex problem, the main element to point-out about this dynamic is its Hamiltonian nature. The Hamiltonian of the point-vortex system is given by
\begin{equation}\label{def:hamiltonnian}
\begin{array}{cccc}
H:&\RR^{2N}&\longrightarrow&\RR,\\\quad &X=(x_1\dots x_N)&\longmapsto&\displaystyle\sum_{i\neq j}a_i\,a_j\, G\big(|x_i-x_j|\big).
\end{array}
\end{equation}
The system~\eqref{eq:vortex equation} can be rewritten 
\begin{equation}\label{eq:hamiltonnian formulation}
a_i\frac{d}{dt}x_i(t)=\nabla^\perp_{x_i}H(X).
\end{equation}
The first consequence of this Hamiltonian reformulation is the preservation of the Hamiltonian $H$ along the flow $S^t$ of~\eqref{eq:vortex equation}:
\begin{equation}\label{lem:preservation Hamiltonien}
\forall\;t\in[0,T),\qquad\frac{d}{dt}H(S^tX)=0.
\end{equation}
We recall that the flow of a differential equation is the function $S^t$ that maps the position $X\in\RR^{2N}$ at time $t=0$ to the position at time $t$. In other words, $S^tX=\big(x_1(t),\dots,x_N(t)\big)$ solution to~\eqref{eq:vortex equation}, with initial positions $X=(x_1,\dots,x_N)$.
Another consequence of the Hamiltonian of the system is the Liouville theorem that ensures the preservation of the Lebesgue measure by the flow.
More precisely, if $V_0\subseteq\RR^{2N}\setminus\mC_T$ is measurable, then we have 
\begin{equation}\label{lem:liouville theorem}
\forall\;t\in[0,T),\qquad\frac{d}{dt}\;\cL^{2N}(S^tV_0)=0,\end{equation}
where $\cL^{2N}$ denotes the Lebesgue measure on $\RR^{2N}$.
For the proof of the Liouville Theorem, we refer to the one given by Arnold in~\cite[Part 3]{Arnold_1978}.
With the Hamiltonian formulation also comes the Noether theorem~\cite{Noether_1918} that provides the quantities left invariant by the flow corresponding to the geometrical invariances  of the Hamiltonian $H$.
The \emph{vorticity vector} is defined for all initial datum $X\in\RR^{2N}$ by
\begin{equation}\label{def:vorticity vector}
M(X):=\sum_{i=1}^Na_i\,x_i.
\end{equation}
The translations invariance of $H$ implies the conservation of the vorticity vector:
\begin{equation}\label{lem:vorticity vector}
\forall\;t\in[0,T),\qquad\frac{d}{dt}M(S^tX)=0.\end{equation}
When the system is non-neutral, meaning that $\sum_ia_i\neq0$, this lemma implies the preservation of the \emph{center of vorticity} of the system defined by
\begin{equation}\label{def:vorticity center}
B(X):=\Big(\sum_{i=1}^Na_i\Big)^{-1}\sum_{i=1}^Na_i\,x_i.
\end{equation}
Similarly, the invariance by the rotations, implies the conservation of the moment of inertia defined  by
\begin{equation}\label{def:inertia momentum}
I(X):=\sum_{i=1}^Na_i\,|x_i|^2.
\end{equation}
We have:
\begin{equation}\label{lem:inertia momentum}
\forall\;t\in[0,T),\qquad\frac{d}{dt}I(S^tX)=0.\end{equation}
The combination of these two lemmas implies the preservation of
\begin{equation}\label{def:collapse constraint}
C(X):=\sum_{i=1}^N\sum_{\substack{j=1\\i\neq j}}^Na_i\,a_j\,|x_i-x_j|^2.
\end{equation}
Indeed, if we expand the square in the right-hand side of~\eqref{def:collapse constraint}, we obtain by a straight-forward calculation
\begin{equation}\label{okay}C(X)=2\Big(\sum_{i=1}^Na_i\Big)I(X)-2\big|M(X)\big|^2.
\end{equation}
The preservation of this quantity is referred as a collapse constraint because it is widely used in the study of vortex collapses for small number of vortices~\cite{Novikov_1975, Novikov_Sedov_1979, Aref_1979, Badin_Barry_2018}. 
Indeed, a collapse means that $|x_i-x_j|^2$ vanishes for some values of $i$ and $j$. 
Combined with the preservation of $C$, this gives a necessary condition for a vortex collapse. For instance, in the case of a collapse for a system of $3$ vortices, this gives the constraint $C=0$.

\section{Main results}

\subsection{Uniform bound results}

\subsubsection{The uniform bound Theorem}
The specific case of the Euler point-vortex system corresponds to the Green function of the Laplacian~\eqref{def:Green functions}.
This particular case is studied in \cite[Chap.\,4]{Marchioro_Pulvirenti_1994}. 
More precisely, they focused on a specific situation for which the intensities for the vortices satisfy 
\begin{equation}\label{eq:no null partial sum}
\forall\;A\subseteq\{1\dots N\}\; s.t.\;\;A\neq\emptyset,\qquad\sum_{i\in A}a_i\neq0.
\end{equation}
A vortex system such that the sum of all the intensities $a_i$ is equal to $0$ is called in~\cite[Chap.\,4]{Marchioro_Pulvirenti_1994} a ``\emph{neutral system}''. 
No name to Hypothesis~\eqref{eq:no null partial sum} is given and we suggest that to call it ``\emph{non-neutral clusters hypothesis}''.

The main interest of this hypothesis relies on the preservation of the center of vorticity property~\eqref{lem:vorticity vector}.
Under the non-neutral cluster hypothesis, the center of vorticity is well defined, not only for the whole system but also for any subset of vortices.
It can be said intuitively that a vortex cluster is expected to ``\emph{turn around its center of vorticity}''.
More precisely, we provide a bound on the trajectories that is uniform with respect to the initial datum $X\in\RR^{2N}$ but also with respect to the singularity of the kernel profile $G$ near $0$. 

\begin{proposition}[Uniform bound on the trajectories]\label{thrm:borne uniforme}
Consider the point-vortex dynamic~\eqref{eq:vortex equation} under the non-neutral clusters hypothesis~\eqref{eq:no null partial sum} with a kernel profile $G\in\cC^{1,1}_{loc}\big(\RR_+^\ast\big)\cap\cC^{1,1}\big([1,+\infty[\big)$.

Then, given any positive time $T>0$, there exist a constant $C$ such that for all initial datum $X\in\RR^{2N}\setminus\mC_T$,
\begin{equation}\label{eq:borne uniforme}
\sup\limits_{t\in[0,T)}\big|X-S_G^tX\big|\leq C,
\end{equation}
where $S^t_G$ is the~\eqref{eq:vortex equation} flow associated to kernel profile $G$. Moreover, the constant $C$ depends only on $N$, the intensities $a_i$, the final time $T$ and on supremum of $r\mapsto\big|\frac{dG}{dr}(r)\big|$ for $r\geq 1$. This constant $C$ does not depend on the initial datum $X\in\RR^{2N}$ nor on the singularity of the kernel $r\mapsto G(r)$ when $r\to0$.
\end{proposition}

In~\cite[Chap.\,4]{Marchioro_Pulvirenti_1994}, a weaker version of this result is established only for the Euler case.
We extend this result to a general case where it holds no matter what the singularity of the kernel $G$ in $0^+$ is with a proof widely inspired from Theorem $2.1$ in~\cite[Chap.\,4]{Marchioro_Pulvirenti_1994}.
We remark that the non-neutral cluster hypothesis~\eqref{eq:no null partial sum} is essential. 
Indeed, the simple situation of a vortex pair with intensities $+1$ and $-1$ gives raise to a translation motion along two parallel lines at a speed that blows up as the initial distance between the two vortices goes to $0^+$. 
The trajectories are bounded in finite time but the bound is not uniform, depending on the initial conditions and on the singularity of the kernel near $r=0$.
We underline that Theorem~\ref{thrm:borne uniforme} apply to the quasi-geostrophic case given by 
the Green function of the fractional Laplacian~\eqref{def:Green functions}. 
In this sense, this theorem is the 
extension of Theorem $2.1$ in~\cite[Chap.\,4]{Marchioro_Pulvirenti_1994} 
to the quasi-geostrophic case. 
Another improvement lays in an explicit computation of the constant appearing in~\eqref{eq:borne uniforme} (unlike~\cite[Chap.\,4]{Marchioro_Pulvirenti_1994}). It is given by a recursive formula with respect to the number off vortices.

\subsubsection{The case of intensities $a_i$ all positive}
As a further consequence of Theorem~\ref{thrm:borne uniforme} we can show the the impossibility of collapses in the case where the intensities $a_i$ are all positive.

\begin{corollary}\label{lem:non collapse}
Let $G$ be a profile such that
\begin{equation}
\big|G(r)\big|\longrightarrow+\infty\qquad\text{as } r\to0.\label{eq:hyp on G 1}
\end{equation}
Assume that the intensities $a_i$ are all positive and consider an
initial datum $X\in\RR^{2N}$ for the point-vortex system~\eqref{eq:vortex equation} such that $x_i\neq x_j$ for all $i\neq j$. Then there is no collapse of vortices at any time.
\end{corollary}
Contrarily to the existence result given by Theorem~$2.2$ 
in~\cite[Chap.\,4]{Marchioro_Pulvirenti_1994} and its generalization to (SQG) at Theorem~\ref{thrm:Improved Marchioro Pulvirenti}, 
this result is true for all initial datum and not only for almost every one. Hypothesis~\eqref{eq:hyp on G 1} may appear a bit restrictive. Indeed, if we consider the kernels
\begin{equation}
\forall\;r>0,\qquad\frac{dG}{dr}(r)=\frac{1}{r^\alpha},
\end{equation} 
with $0<\alpha<1$ then the associated kernel $G$ does not satisfy~\eqref{eq:hyp on G 1}. 
The possibility to extend Corollary~\ref{lem:non collapse} to this case is an open problem. 
Nevertheless, the physical relevant cases are $\alpha=1$ for the Euler model or $3<\alpha<1$ for the quasi-geostrophic model. 
For these values of $\alpha$, Corollary~\ref{lem:non collapse} apply.

\subsubsection{The uniform relative bound theorem}
A natural question concerning Theorem~\ref{thrm:borne uniforme} is to ask what this result becomes when the non-neutral cluster hypothesis~\eqref{eq:no null partial sum} ceases to be satisfied. For instance, 
we consider instead of~\eqref{eq:no null partial sum} the following hypothesis:
\begin{equation}\label{eq:no null sub partial sum}
\forall\;A\subseteq\{1\dots N\}\; s.t.\;\;A\neq\emptyset\;\;\text{or}\;\;\{1\dots N\},\qquad\quad\sum_{i\in A}a_i\neq0.
\end{equation}
In other words, all the strict sub-clusters must have the sum of their 
intensities different from $0$ but we allow the 
total sum $\sum_{i=1}^Na_i$ to be equal to $0$. This situation is 
achieved for instance by the vortex pair of intensities $+1$ 
and $-1$ that are translating at a constant speed. 

\begin{proposition}[Uniform relative bound on the trajectories]\label{thrm:uniform relative bound}
For a given set of points noted $X=(x_1\dots x_N)\in\RR^{2N}$, we define the diameter of this set by
\begin{equation}\label{def:diam}
diam(X)\;:=\;\max\limits_{i\neq j}|x_i-x_j|.
\end{equation}

Consider the point-vortex dynamic~\eqref{eq:vortex equation} under hypothesis~\eqref{eq:no null sub partial sum} with a kernel profile $G\in\cC^{1,1}_{loc}\big(\RR_+^\ast\big)\cap\cC^{1,1}\big([1,+\infty[\big)$.
Let $T>0$ the final time.
Then for all kernel profile $G\in\cC^{1,1}_{loc}\big(\RR_+^\ast\big)\cap\cC^{1,1}\big([1,+\infty[\big)$ and for all
initial datum $X\in\RR^{2N}$ that are not leading to collapse on $[0,T)$,
\begin{equation}\label{eq:borne uniforme relative}
\sup\limits_{t\in[0,T)}\diam\big(S_G^tX\big)\leq \diam(X)+C,
\end{equation}
where $S^t_G$ is the flow associated to~\eqref{eq:vortex equation} 
with the kernel profile equal to $G$.
Moreover, the constant $C$ depends only on $N$, the intensities $a_i$, the final time $T$ and on supremum of $r\mapsto\big|\frac{dG}{dr}(r)\big|$ for $r\geq 1$. 
\end{proposition}

This constant $C$ does not depend on the initial datum $X\in\RR^{2N}$ nor on the singularity of the kernel $r\mapsto G(r)$ when $r\to0$. The reasoning is close to the one of Proposition~\ref{thrm:borne uniforme}.

Similarly, the constant is computed explicitly  throughout the proof of the theorem.

\subsection{Convergence result for Euler point-vortices}\label{sec:convergence result}
The systems of vortices that are studied with more details 
in~\cite[Chap.\,4]{Marchioro_Pulvirenti_1994} are the systems for which 
the non-neutral clusters hypothesis~\eqref{eq:no null partial sum} 
holds. As stated in the previous section, their result concerning 
uniform bound on the trajectories can be improved to consider a much 
wider class of point-vortex systems, including the quasi-geostrophic 
case. Nevertheless, in the particular case of the Euler point-vortex 
system with the non-neutral clusters hypothesis~\eqref{eq:no null partial sum}, we are also able to write a convergence result.
When vortices come to collapse, their speed may become infinite as a 
consequence of the kernel profile singularity in $0^+$. 
But if their speed blows up, then any pathological behavior near the 
time of collapse $T_X$ is \emph{a priori} possible. 
We prove here that the trajectories are actually convergent in the Euler case.

\begin{theorem}[Convergence for Euler vortices under non-neutral clusters hypothesis]\label{thrm:generalization of no partial sum with ordinary}
Consider the point-vortex model~\eqref{eq:vortex equation} under 
hypothesis~\eqref{eq:no null partial sum} with a kernel profile $G_1$ 
corresponding to the Green function of the 
Laplacian~\eqref{def:Green functions}.
Let $X\in\mC$ be an initial datum leading to a collapse at time $T_X$. 
Then, for all $i=1\dots N$, there exists an $x_i^\ast\in\RR^{2}$ such that
\begin{equation}\label{eq:continuity of the trajectories 2}
x_i(t)\longrightarrow x_i^\ast\qquad\text{as }t\to T_X^-.
\end{equation}
\end{theorem}
The first step of the proof consists in the following idea:
For a fixed value of $t\in[0,T_X)$ define the distribution
$P_t:=\sum_{i=1}^Na_i\delta_{x_i(t)}.$
The point-vortex equation gives that this distribution converges as $t\to T_X$.
In a second time, it is possible to prove that this convergence is 
actually stronger and obtain~\eqref{eq:continuity of the 
trajectories 2} by exploiting the non-neutral cluster hypothesis.
The proof of this result is specific to the Euler case and we do not 
know yet whether the conclusion extends to the (SQG) case.

\subsection{Improbability of collapses for point-vortices}\label{sec:improbability collapses}
Understanding the sets of collapses $\mC_T$ is an important issue for the study of point-vortices since these sets give the time of existence of the point-vortex system for a given initial datum.
Marchioro and Pulvirenti in their study of the point-vortex 
problem~\cite{Marchioro_Pulvirenti_1994} provided the following 
improbability result

\begin{theorem}[Improbability of collapses, Marchioro-Pulvirenti, 1993]~\label{thrm:Marchioro Pulvirenti existence theorem}
Consider the point-vortex problem~\eqref{eq:vortex equation} with 
kernel profile $G_1$, the kernel associated to the Green function of 
the Laplacian on the plane~\eqref{def:Green functions}. Assume that the intensities of the vortices satisfy
the \emph{non neutral cluster hypothesis}~\eqref{eq:no null partial sum}.

Then, the set of initial datum for the dynamic~\eqref{eq:vortex equation} that lead
to collapses in finite time
is a set of Lebesgue measure equal to $0$.
\end{theorem}
The non neutral cluster hypothesis~\eqref{eq:no null partial sum} is an important hypothesis for this theorem as it allows us to use Theorem~\ref{thrm:borne uniforme} that provides a uniform bound.
Concerning this result, the most natural question consists in 
removing this assumption on the intensities of the vortices, which eventually leads to the 
following conjecture. 
\begin{conjecture}\label{conj:improbability of collapses}
The set of collapses $\mC$ has a Lebesgue measure $0$.
\end{conjecture}
It is not yet possible to prove such a conjecture.
The main difficulty lays in the understanding of the situations 
where some vortices collide in such a way that they go to infinity 
in finite time, or show an unbounded pathological behavior. The existence of unbounded trajectories in finite time is also an open problem.
Although we are not able to answer to Conjecture~\eqref{conj:improbability of collapses}, 
we are able to improve Theorem~\ref{thrm:Marchioro Pulvirenti existence theorem} as stated in the following theorem:
\begin{theorem}[improbability of collapses for Euler and SQG vortices]\label{thrm:Improved Marchioro Pulvirenti}
Consider the point-vortex problem~\eqref{eq:vortex equation} with 
kernel profile $G_s$, the kernel associated to the Green function of 
the Laplacian or the fractional Laplacian on the plane~\eqref{def:Green functions} for $s\in(0,1]$. Assume that the intensities of the vortices satisfy~\eqref{eq:no null sub partial sum}.

Then, the set of initial datum leading to collapses 
has a Lebesgue measure equal to $0$.
\end{theorem}

This theorem is slightly more general than 
Theorem~\ref{thrm:Marchioro Pulvirenti existence theorem} for two reasons. 
First, it is true both for Euler and for quasi-geostrophic point-vortex 
models. This aspect was already partially improved by~\cite{Geldhauser_Romito_2020}, where they obtained the result
for $s>1/2$. Indeed, the value $s=1/2$ appear to be a critical value for the quasi-geostrophic equations where 
integrability problems arise. Our arguments manage to pass trough these difficulties and obtain 
the result for all $s\in(0,1]$.
The second improvement lays in the fact that we managed to replace the non neutral cluster 
hypothesis~\eqref{eq:no null partial sum} by the weaker hypothesis~\eqref{eq:no null sub partial sum}. 
This weaker hypothesis can be seen at first sight as a small improvement. 
Nevertheless, it make a quite important difference because the non neutral cluster 
hypothesis~\eqref{eq:no null partial sum} implies that the trajectories are bounded (Theorem~\ref{thrm:borne uniforme}) 
whereas hypothesis~\eqref{eq:no null sub partial sum} only imply a relative bound (Theorem~\ref{thrm:uniform relative bound}).
In other words, this improbability result allows a simple unbounded behaviors: the cases where the vortices collectively goes to infinity.
In the article~\cite{Geldhauser_Romito_2020}, the authors study the evolution in time of the following function:
\begin{equation}
\Phi(X):=\sum_{i=1}^N\sum_{\substack{j=1\\j\neq i}}^N\big|x_i-x_j\big|^{-\beta}
\end{equation}
when the $x_i(t)$ evolves according to the point-vortex equation~\eqref{eq:vortex equation} to obtain sufficent conditions for a collapse. 
Indeed, this function blows up if $|x_i(t)-x_j(t)|\to0$ as $t\to T$ at a speed depending on the choice of $\beta>0$ and then it is possible to study the collapses as being the set points for wich this function becomes unbounded in time.
Unfortunately, the choice of $\beta$ that must be made for this proof (depending on the value of $s$) create intergability problems when $0<s\leq1/2$.
We overcome these problems arising in the proof of~\cite{Geldhauser_Romito_2020} when $s\leq1/2$ by replacing the singularity $r^{-\beta}$ in the definition of $\Phi$ by a regularized version with parameter $\varepsilon>0$. 
We then proceed to a reasoning similar to the one of Marchioro and Pulvirenti~\cite[Chap.\,4]{Marchioro_Pulvirenti_1994} and then conclude by letting $\varepsilon\to0.$

\section{Outlines of the proofs}
To ease the general reading and understanding of the article, we 
draw here only the outlines of the proofs of the main Theorems. These larger 
proofs are decomposed into smaller lemmas stated in this section and 
forming the main intermediate steps. 
The links and articulations between the different lemmas are
developed in this section but the detailed technical proofs of each 
different Lemma are all postponed to Section~\ref{sec:proofs}. 
Concerning Propositions~\ref{thrm:borne uniforme} and~\ref{thrm:uniform relative bound}, 
the proofs are shorter and directly done in section~\ref{sec:proofs}.

\subsection{Outline of the proof for Theorem~\ref{thrm:generalization of no partial sum with ordinary}}
The first part of the proof consists in considering Dirac masses 
located on the vortices and to prove that this converges in the 
distributional sense as $t\to T_X$, the time of collapse.

\begin{lemma}[Convergence of the Dirac measures]\label{lem:Dirac}
Let $X\in\RR^{2N}$. From the vortices $x_i(t)$ evolving according to 
the differential equations~\eqref{eq:vortex equation} with initial 
datum $X$, define the distribution
\begin{equation}\label{def:Dirac sum}
P_X(t):=\sum_{i=1}^Na_i\,\delta_{x_i(t)},
\end{equation}
where $\delta_x$ denotes the Dirac mass at point $x\in\RR^2$. 
Assume now that the evolution problem associated to the initial 
datum $X$ is well defined on an interval of time $[0,T)$ for some $T>0$. 
Then, there exists $X^\ast\in\RR^{2N}$ and $b\in\{0,1\}^N$ such that
\begin{equation}\label{eq:Dirac limit}
\sum_{i=1}^Na_i\,\delta_{x_i(t)}\longrightarrow\sum_{i=1}^Na_i\,b_i
\,\delta_{x_i^\ast}\qquad\text{in the weak sense of measure as } t\to T^-
\end{equation}
and such that  $$b_i=0\qquad\Longrightarrow
\qquad\sup\limits_{t\in[0,T[}|x_i(t)|=+\infty.$$
\end{lemma}

Note that the proof makes no use of the non-neutral clusters 
hypothesis~\eqref{eq:no null partial sum}. In the particular case 
where this hypothesis is satisfied, theorem~\ref{thrm:borne uniforme} implies that the vortices stay bounded on bounded 
intervals of time. Therefore, in this particular case all the 
coefficients $b_i$ given by this lemma are equal to $1$.

The next step of the proof exploits more precisely Hypothesis~\eqref{eq:no null partial sum} and the continuity of the 
trajectories to obtain that a given vortex $x_i(t)$ can only 
converge as $t\to T^-$ to one element of the set $\{x_1^\ast\dots x_N^\ast\}$ given by Lemma~\ref{lem:Dirac}. More precisely, if a 
given point $x_i(t)$ have at least two adherence points, it is 
possible to extract a third adherence point that does not belongs to 
$\{x_1^\ast\dots x_N^\ast\}$ and this provides a contradiction with 
Lemma~\ref{lem:Dirac}. See Section~\ref{sec:proofs} for a detailed proof.

\subsection{Outline of the proof for Theorem~\ref{thrm:Improved Marchioro Pulvirenti}}

\subsubsection{The modified system}
Let $i$ be fixed in $\{1\dots N\}.$
The modified system consists in studying the evolution of $y_{ij}:=x_i-x_j$ for $j\in\{1\dots N\}\setminus\{i\}$. The idea is that knowing the relative position of the vortices (the differences $y_{ij}$) is enough to study the problem of collapses.
The initial problem~\eqref{eq:vortex equation} implies that
\begin{equation}\label{eq:evolution difference}
\frac{d}{dt}(x_i-x_j)(t)=\sum_{k\neq i}a_k\nabla^\perp G_s\big(|x_i-x_k|\big)-\sum_{k\neq j}a_k\nabla^\perp G_s\big(|x_j-x_k|\big).
\end{equation}
Therefore, the evolution of $y_{ij}$ is given by
\begin{equation}
\label{eq:evolution Y}
\frac{d}{dt}y_{ij}=(a_i+a_j)\nabla^\perp G_s\big(|y_{ij}|\big)+
\sum_{k\neq i,j}a_k\Big(\nabla^\perp G_s\big(|y_{ik}|\big)
+\nabla^\perp G_s\big(|y_{ij}-y_{ik}|\big)\Big).
\end{equation}
The main interest of this new system is that Theorem~\ref{thrm:Improved Marchioro Pulvirenti} can be reformulated using 
only the differences $y_{ij}$. 

\begin{lemma}[Reformulation of Theorem~\ref{thrm:Improved Marchioro Pulvirenti}]\label{lem:reformulation}
Denote by $Y_i(t):=\big(y_{ij}(t)\big)_{j\neq i}$ the solution 
to~\eqref{eq:evolution Y} at time $t$ with initial datum $Y_i\in\RR^{2(N-1)}$.
Assume that for all $i\in\{1\dots N\}$ and for all $T>0$ and $\rho>0$ the set
\begin{equation}\label{eq:the set with measure zero}
\Big\{Y_i:=(y_{ij})_{j\neq i}\in\cB(0,\rho)^{2(N-1)}:\exists\;T_X\in[0,T],
\quad\liminf_{t\to T_X^-}\;\min\limits_{j\neq i}\big|y_{ij}(t)\big|=0\Big\}.
\end{equation}
has its Lebesgue measure $\cL^{2(N-1)}$ equal to $0$. Then in this 
case the conclusion of Theorem~\ref{thrm:Improved Marchioro Pulvirenti} holds.
\end{lemma}
The notation $\cB(x_0,\rho)$ refers to the Euclidean ball of $\RR^2$ of center $x_0$ and radius $\rho$.
The rest of the work consists then in studying this 
system~\eqref{eq:evolution Y} and to establish 
that~\eqref{eq:the set with measure zero} does have measure $0$.
This modified dynamics has many properties in common with the 
original one. In particular, this 
new dynamics still satisfies the Liouville property. Define the 
function $\cH_{ij}:\RR^{2(N-1)}\to\RR$ by
\begin{equation}\label{eq:hamilton Y}
\cH_{ij}\Big[(y_{ik})_{k\neq i}\Big]:=(a_i+a_j)G_s\big(|y_{ij}|\big)
+y_{ij}\cdot\!\sum_{k\neq i,j}a_k\nabla G_s\big(|y_{ik}|\big)+
\sum_{k\neq i,j}a_kG\big(|y_{ij}-y_{ik}|\big).
\end{equation}
Combining this equation with~\eqref{eq:evolution Y} gives
\begin{equation}\label{eq:nabla perp Y}
\frac{d}{dt}y_{ij}=\nabla^\perp_{y_{ij}}\cH_{ij}\Big[(y_{ik})_{k\neq i}\Big].
\end{equation}
This equation says that, in a certain sense, the dynamic of the 
vector $Y_i:=(y_{ij})_{j\neq i}$ shows an Hamiltonian structure. It 
is not an Hamiltonian system because the function $\cH_{ij}$ does 
depend on $j$ but there is still a structure with an operator $\nabla^\perp$. Therefore, the Schwartz theorem gives
\begin{equation}\label{eq:Schwarz Y}
\di_{y_{ij}}\Big(\frac{d}{dt}y_{ij}\Big)=\di_{y_{ij}}\bigg(\nabla^\perp_{y_{ij}}\cH_{ij}\Big[(y_{ik})_{k\neq i}\Big]\bigg)=0.
\end{equation}
Since the velocity is divergent-free, a Liouville theorem holds for this dynamics. 

\begin{lemma}[Liouville theorem]\label{lem:Liouville Theorem for the modified dynamics}
The Lebesgue measure on the space $\RR^{2(N-1)}$ given by
\begin{equation}
\prod_{j\neq i}dy_{ij}
\end{equation}
is preserved by the flow $\mS_{i}^t$ associated to~\eqref{eq:nabla perp Y}.
\end{lemma}

A detailed proof of the Liouville theorem in a more general setting 
can be found at~\cite[Part 3]{Arnold_1978} and we refer to it for the 
proof of Lemma~\ref{lem:Liouville Theorem for the modified dynamics}. 

\subsubsection{Estimate the collapses} 
From the kernel $G_s$ given at~\eqref{def:Green functions}, define 
now the regularized profiles $G_{s,\varepsilon}$ for 
$\varepsilon\in(0,1]$. The objective here is to drop the singularity of 
the kernel $G_s$ near $0^+$. We ask $G_{s,\varepsilon}$ to be 
$\cC^1$ on $\RR_+$ (until the boundary) and to verify
\begin{align}
&\bullet\qquad G_{s,\varepsilon}(q)=G_s(q)\qquad \mathrm{when}\quad \varepsilon\leq q,\label{eq:condition on G epsilon 1}\\
&\bullet\qquad |G_{s,\varepsilon}(q)|\leq |G_s(q)|\qquad\mathrm{for\;all}\quad q\in\RR_+,\label{eq:condition on G epsilon 2}\\
&\bullet\qquad \Big|\frac{d}{dq}G_{s,\varepsilon}(q)\Big|
\leq \Big|\frac{d}{dq}G_s(\varepsilon)\Big|\qquad\mathrm{for\;all}\quad q\leq\varepsilon,\label{eq:condition on G epsilon 3}\\
&\bullet\qquad |G_{s,\varepsilon}(q)|\leq 2|G_s(\varepsilon)|\qquad\mathrm{for\;all}\quad q\in\RR_+.\label{eq:condition on G epsilon 4}
\end{align}
Since $G_{s,\varepsilon}$ is of class $\cC^{1,1}$ on $\RR_+$, the 
dynamic defined by~\eqref{eq:vortex equation} for the kernel 
$G_{s,\varepsilon}$ is always well-defined and is Hamiltonian.
In the sequel we denote $\mS_{i,\varepsilon}^t$ the flow at time $t$ 
associated to the evolution equation~\eqref{eq:evolution Y} with the 
kernel profile $G_s$ replaced by the regularization $G_{s,\varepsilon}$. 
The motion induced by the kernel profile $G_{s,\varepsilon}$ 
coincides with the original motion provided that the distances 
between the vortices remain higher than $\varepsilon$. 
Theorem~\ref{thrm:Improved Marchioro Pulvirenti} can be reformulated 
as follows.

\begin{lemma}[Reformulation of 
Theorem~\ref{thrm:Improved Marchioro Pulvirenti} with $\varepsilon$-Regularized dynamic] \label{lem:reformulation convergence}
Assume that for all $T>0$ and for all $\rho>0$ we have the following convergence:
\begin{equation}\label{eq:reformulation}
\cL^{2(N-1)}\Big\{Y_i=(y_{ij})_{j\neq i}\in\cB(0,\rho)^{2(N-1)}\;:
\;\min\limits_{j\neq i}
\inf\limits_{t\in[0,T]}\big|y_{ij}^\varepsilon(t)\big|\leq\varepsilon\Big\}
\xrightarrow[\varepsilon\to0^+]{}0.
\end{equation}
Then for all $i\in\{1\dots N\}$, the set~\eqref{eq:the set with measure zero} has Lebesgue measure $0$.
\end{lemma}

This lemma is nothing more than Lemma~\ref{lem:reformulation} where 
we added this parameter $\varepsilon>0$ allowing us to regularize 
the kernel.
Therefore, combined with Lemma~\ref{lem:reformulation}, this lemma 
gives Theorem~\ref{thrm:Improved Marchioro Pulvirenti} provided that 
the convergence~\eqref{eq:reformulation} holds. 

\begin{lemma}\label{lem:collapses}
Let $i\in\{1\dots N\}$, $\varepsilon>0$ and $\rho>0$. Then
\begin{equation}\label{eq:convergence}
\cL^{2(N-1)}\Big\{Y_i=(y_{ij})_{j\neq i}\in\cB(0,\rho)^{2(N-1)}\;:
\;\min\limits_{j\neq i}\inf_{t\in[0,T]}\big|y_{ij}^\varepsilon(t)\big|\leq\varepsilon\Big\}
\leq C\left\{\begin{array}{ll}
\varepsilon&\quad\text{if }s>0.5,\\
\varepsilon\log(1/\varepsilon)&\quad\text{if } s=0.5,\\
\varepsilon^{2s}&\quad\text{if } s<0.5.
\end{array}\right.
\end{equation}
where the constant $C$ only depends on $N$, $|a_i|$, $T$ and $\rho$.
\end{lemma}
The proof of this last lemma is reminiscent from the proof of 
Theorem~\ref{thrm:Marchioro Pulvirenti existence theorem} with some 
new arguments. The main idea is to rely on a Bienaymé-Tchebycheff 
inequality applied to a well-chosen function. 
This last estimate~\eqref{eq:convergence} with Lemma~\ref{lem:reformulation convergence} 
concludes the proof of Theorem~\ref{thrm:Improved Marchioro Pulvirenti}.\qed

\section{Technical proofs}\label{sec:proofs}

\subsection{Proofs for the Hamiltonian formulation}

\subsubsection{Proof of the preservation of the Hamiltonian~\eqref{lem:preservation Hamiltonien}}
Let $X\in\RR^{2N}\setminus\mC_T$. The Hamiltonian formulation~\eqref{eq:hamiltonnian formulation} implies
\begin{equation}
\frac{d}{dt}H(S^tX)=\sum_{i=1}^N\nabla_{x_i}H(S^tX)\cdot\frac{d}{dt}x_i(t)=\sum_{i=1}^N\frac{1}{a_i}\nabla_{x_i}H(S^tX)\cdot\nabla_{x_i}^\perp H(S^tX)=0,
\end{equation}
where we used the general fact: $\nabla_xf(x)\cdot\nabla^\perp_xf(x)=0.$
\qed

\subsubsection{Proof of the conservation of the vorticity vector~\eqref{lem:vorticity vector}}
Let $X\in\RR^{2N}\setminus\mC_T$, by direct computation using~\eqref{eq:vortex equation}:
\begin{equation}
\frac{d}{dt}M(S^tX)=\sum_{i=1}^Na_i\,\frac{d}{dt}x_i(t)=\sum_{i=1}^N\sum_{\substack{j=1\\j\neq i}}a_i\,a_j\,\nabla^\perp G\big(|x_i(t)-x_j(t)|\big)=0,
\end{equation}
where we used the identity $\nabla^\perp G(|-x|)=-\nabla^\perp G(|x|).$ \qed

\subsubsection{Proof of the conservation of the moment of inertia~\eqref{lem:inertia momentum}}
Let $X\in\RR^{2N}\setminus\mC_T$, by direct computation using~\eqref{eq:vortex equation}:
\begin{equation}
\frac{d}{dt}I(S^tX)=2\sum_{i=1}^Na_i\,x_i(t)\cdot\frac{d}{dt}x_i(t)=2\sum_{i=1}^N\sum_{\substack{j=1\\j\neq i}}^Na_i\,a_j\,x_i(t)\cdot\nabla^\perp G\big(|x_i(t)-x_j(t)|\big).
\end{equation}
If we proceed to a symmetrization of the double sum above, swapping the indices $i\leftrightarrow j$, 
and using the identity $\nabla^\perp G(|-x|)=-\nabla^\perp G(|x|)$, we are led to
\begin{equation}
\frac{d}{dt}I(S^tX)=\sum_{i=1}^N\sum_{\substack{j=1\\j\neq i}}^Na_i\,a_j\,\big(x_i(t)-x_j(t)\big)\cdot\nabla^\perp G\big(|x_i(t)-x_j(t)|\big).
\end{equation}
Observing now that the vector $\nabla^\perp G(|x|)$ is orthogonal to the vector $x$, we deduce that all the scalar products appearing in the expression above are $0$.\qed

\subsection{Proofs for the uniform bound results}

\subsubsection{Proof of Propositions~\ref{thrm:borne uniforme} and~\ref{thrm:uniform relative bound}}
The two uniform bounds given by Propositions~\ref{thrm:borne uniforme} and~\ref{thrm:uniform relative bound} are a consequence of the following proposition:

\begin{proposition}\label{prop:borne reformule}
Let $N\in\NN$ with $N\neq0$ and let $a_i\neq0$ with $i=1\dots N$. let $X\in\RR^{2N}$ be an initial datum for the following differential evolution problem for $t\in[0,T)$:
\begin{equation}\label{eq:vortex perturbed}
\frac{d}{dt}x_i(t)\;=\;\sum_{\substack{j=1\\j\neq i}}^Na_j\,\nabla^\perp G\big(|x_i(t)-x_j(t)|\big)+f_i\big(t,x_i(t)\big),
\end{equation}
where $G\in\cC^{1,1}_{loc}\big(\RR_+^\ast\big)\cap\cC^{1,1}\big([1,+\infty[\big)$ 
is the kernel profile and where $f_i:[0,T)\times\RR^2\to\RR^2$ are smooth external fields. 
One makes the assumption that there are no collapses for all $t\in[0,T)$ so that the dynamic~\eqref{eq:vortex perturbed} is well-defined.
One sets:
\begin{equation}\label{def:a A0 A}\begin{split}
&a:=\sum_{i=1}^N|a_i|,\\
&A_0:=\min\limits_{\substack{\cP\subseteq\{1\dots N\}\\\cP\neq\emptyset,\,\cP\neq\{1\dots N\}}}\bigg|\sum_{i\in\cP}a_i\bigg|,\\
&A:=\min\limits_{\substack{\cP\subseteq\{1\dots N\}\\\cP\neq\emptyset}}\bigg|\sum_{i\in\cP}a_i\bigg|\;=\min\Big\{A_0;\Big|\sum_{i=1}^Na_i\Big|\Big\}.
\end{split}
\end{equation}

Then there exists a function $C:\NN\times(\RR_+)^5\to\RR_+$, given explicitly, that is non-decreasing with respect to any of its $6$ variables and such that:\vspace{0.2cm}

$(i)$ If $A_0\neq0$, then $\forall\;t\in[0,T),\;\forall\, i,j\in\{1\dots N\},$
\begin{equation}\label{eq:announced estimate}
\Big|\big(x_i(t)-x_j(t)\big)-\big(x_i(0)-x_j(0)\big)\Big|\;\leq\;C\Big(N,\,\max_k\|f_k\|_{L^\infty},\,T,\,a,\,\frac{1}{A_0},\,\sup_{r\geq1}\bigg|\frac{dG}{dr}(r)\bigg|\Big).
\end{equation}

$(ii)$ If moreover $A\neq 0$, then $\forall\;t\in[0,T),\;\forall\, i\in\{1\dots N\},$
\begin{equation}
\Big|\big(x_i(t)-B(t)\big)-\big(x_i(0)-B(0)\big)\Big|\;\leq\;\frac{a}{A}\,C\Big(N,\,\max_k\|f_k\|_{L^\infty},\,T,\,a,\,\frac{1}{A_0},\,\sup_{r\geq1}\bigg|\frac{dG}{dr}(r)\bigg|\Big),
\end{equation}
where $B(t)$ is the center of vorticity of the system~\eqref{def:vorticity center}.
\end{proposition}

The first point of this proposition implies the relative uniform bound stated by Proposition~\ref{thrm:uniform relative bound} because the non-neutral sub-clusters hypothesis~\eqref{eq:no null sub partial sum} is equivalent to $A_0\neq0$.
Proposition~\ref{thrm:uniform relative bound} is then obtained by choosing $f_i\equiv0$ for all index $i$.
Similarly, the second point of this proposition implies the uniform bound stated by Proposition~\ref{thrm:borne uniforme} since the non-neutral clusters hypothesis~\eqref{eq:no null partial sum} is equivalent to $A\neq0$. Recall that in the case $f_i\equiv0$, the center of vorticity $B$ is a constant of the movement~\eqref{lem:vorticity vector}.

\begin{proof}
To start with, the second point of this proposition is a direct consequence of the first one because of the following estimate that holds for all $i=1\dots N$:
\begin{equation}
\big|x_i(t)-B(t)\big|=\bigg|\frac{\sum_{j=1}^Na_j\big(x_i(t)-x_j(t)\big)}{\sum_{j=1}^Na_j}\bigg|\leq\frac{a}{A}\,\max_{j,k=1\dots N}|x_j(t)-x_k(t)|.
\end{equation}
Then, there remain to prove Proposition~\ref{prop:borne reformule}-$(i)$. 
The function $C$ is constructed using an iterative argument on the number of vortices $N>0$ that is similar to the proof of the uniform bound in the book of Marchioro and Pulvirenti~\cite[Chap.\,4]{Marchioro_Pulvirenti_1994}.
For the iterative arguments, the case $N=1$ is straight-forward and gives $C(1,\dots)\equiv0$.
Suppose now that the function $C(k,\dots)$ has been constructed for all $k=1\dots N-1$ with $N\geq2$ and satisfy the announced estimate~\eqref{eq:announced estimate}. One is now constructing the function $C(N,\dots)$. For that purpose, one defines in the view of~\eqref{def:a A0 A} the following quantities for all $\cP\subseteq\{1\dots N\}$ non empty:
\begin{equation}
\begin{split}
&a(\cP):=\sum_{i\in\cP}|a_i|,\\
&A_0(\cP):=\min\limits_{\substack{\cQ\subseteq\cP\\\cQ\neq\emptyset,\,\cQ\neq\cP}}\bigg|\sum_{i\in\cQ}a_i\bigg|,\\
&A(\cP):=\min\limits_{\substack{\cQ\subseteq\cP\\\cQ\neq\emptyset}}\bigg|\sum_{i\in\cQ}a_i\bigg|\;=\min\Big\{A_0(\cP);\Big|\sum_{i\in\cP}a_i\Big|\Big\}.
\end{split}
\end{equation}
The following increasing properties hold:
\begin{equation}\label{eq:aA increasing}\begin{split}
&\cQ\subsetneq\cP\qquad\Longrightarrow\qquad a(\cQ)<a(\cP),\\
&\cQ\subsetneq\cP\qquad\Longrightarrow\qquad A(\cQ)\geq A_0(\cP).
\end{split}
\end{equation}
There is also,
\begin{equation}\label{eq:aA and A_0}
\forall\,\cP\subseteq\{1\dots N\},\qquad A_0(\cP)\geq A(\cP).
\end{equation}
Now, let $S>0$ be a parameter supposed very large and that is fixed later on.
One defines the set
\begin{equation}
\mD_S\;:=\;\big\{t\in[0,T):\max_{j,k=1\dots N}|x_j(t)-x_k(t)|\geq S\big\},
\end{equation}
and $t_0:=min\,\mD_S\in[0,T]$. Since the largest distance between vortices is larger or equal to $S$ at time $t_0$, then by the triangular inequality the vortices are divided into two nonempty clusters $\cP,\cQ\subsetneq\{1\dots N\}$ (with $\cP\cup\cQ=\{1\dots N\}$ and $\cP\cap\cQ=\emptyset$) separated by a distance $d$ that is bounded from below by $S/(N-1).$
Indeed, the least favorable case consists in the situation where all 
the vortices forms a rectilinear chain made of $N$ 
points spaced with intervals of same length and that link the two vortices that realize the maximal distance.
If the distance $d$ is larger than $1$ then the interaction between two vortices that does not belong to the same cluster is bounded by $sup_{r\geq1}\big|\frac{dG}{dr}(r)\big|$. 
For that purpose, one makes now the assumption that $S\geq N$ so that $d>1$ (this constraint is handled later in the choice of $S$) n that case, the evolution of the points $x_i$ in cluster $\cP$, using~\eqref{eq:vortex perturbed}, is given by
\begin{equation}\label{eq:vortex cluster evo}
\frac{d}{dt}x_i(t)\;=\;\sum_{\substack{j\in\cP\\j\neq i}}a_j\,\nabla^\perp G\big(|x_i(t)-x_j(t)|\big)+\widetilde{f}_i\big(t,x_i(t)\big),
\end{equation}
where,
\begin{equation}
\widetilde{f}_i\big(t,x_i(t)\big)=f_i\big(t,x_i(t)\big)+\sum_{\substack{j\in\cQ\\j\neq i}}a_j\,\nabla^\perp G\big(|x_i(t)-x_j(t)|\big).
\end{equation}
In particular,
\begin{equation}\label{eq:tilde f}
\|\widetilde{f}_i\|_{L^\infty}\leq\|f_i\|_{L^\infty}+a(\cQ)\,\overline{G'},
\end{equation}
where we have set the notation:
\begin{equation}
\overline{G'}:=\sup_{r\geq1}\bigg|\frac{dG}{dr}(r)\bigg|.
\end{equation}
Then, as long as the distance between the two clusters remain larger than $1$, it is possible to apply the result of Proposition~\ref{prop:borne reformule} recursively to~\eqref{eq:vortex cluster evo}. Note that $A(\cP)>0$ since $A_0=A_0(\{1\dots N\})>0$ by hypothesis and $\cP\neq\{1,\dots,N\}$. It is therefore possible to define the center of vorticity of the cluster $\cP$:
\begin{equation}
B_\cP(t):=\Big(\sum_{j\in\cP}a_j\Big)^{-1}\sum_{j\in\cP}a_j\,x_j(t).
\end{equation}
Proposition~\ref{prop:borne reformule}-$(ii)$ applied recursively to $\cP$ with a number of vortices $\#\cP<N$ gives that for all index $i\in\cP,$
\begin{equation}\label{luth}\begin{split}
\forall\;t\in[t_0,t_1),\qquad\Big|\big(x_i(t)&-B_\cP(t)\big)-\big(x_i(t_0)-B_\cP(t_0)\big)\Big|\\\;&\leq\;\frac{a(\cP)}{A(\cP)}\,C\Big(\#\cP,\,\max_{k\in\cP}\|\widetilde{f}_k\|_{L^\infty},\,|t_1-t_0|,\,a(\cP),\,\frac{1}{A_0(\cP)},\,\overline{G'}\Big),\end{split}
\end{equation}
where
\begin{equation}
t_1:=\sup\big\{t\in[t_0,T)\;:\;\min_{i\in\cP}\,\min_{j\in\cQ}\,|x_i(t)-x_j(t)|\geq1\big\}.
\end{equation}
Note also that $t_1>t_0$ since $S\geq N$ implies that at time $t_0$ the distance is larger than $N/N-1>1$.
Since the function $C$ is increasing with respect to any of its variables, the estimate~\eqref{luth} with the increasing properties~\eqref{eq:aA increasing} becomes, for all $i,j\in\cP,$ and for all $t\in[t_0,t_1)$,
\begin{equation}\label{peter pan}\begin{split}
&\qquad\Big|\big(x_i(t)-B_\cP(t)\big)-\big(x_i(t_0)-B_\cP(t_0)\big)\Big|\qquad\\&\leq\;\frac{a}{A_0}C\Big(N-1,\,\max_k\|\widetilde{f}_k\|_{L^\infty},\,T,\,a,\,\frac{1}{A_0},\,\overline{G'}\Big),\\
&\leq\;\frac{a}{A_0}C\Big(N-1,\,\max_k\|f_k\|_{L^\infty}+a\,\overline{G'},\,T,\,a,\,\frac{1}{A_0},\,\overline{G'}\Big),
\end{split}\end{equation}
Where~\eqref{eq:tilde f} is used for the second inequality.
On the other hand, the center of vorticity is preserved by the flow of the point-vortex dynamics~\eqref{lem:vorticity vector} and therefore $B_\cP$ is only moved by the external field:
\begin{equation}
\forall\;t\in[t_0,t_1),\qquad\frac{d}{dt}B_\cP(t)=\Big(\sum_{j\in\cP}a_j\Big)^{-1}\sum_{j\in\cP}a_j\widetilde{f}\big(t,x_j(t)\big).
\end{equation}
Thus,
\begin{equation}\label{bambi}
\forall\;t\in[t_0,t_1),\qquad\big|B_\cP(t)-B_\cP(t_0)\big|\leq\frac{a(\cP)}{A(\cP)}T\max_k\|\widetilde{f}_k\|_{L^\infty}\leq\frac{a}{A_0}T\Big(\max_k\|f_k\|_{L^\infty}+a\,\overline{G'}\Big),
\end{equation}
where was used~\eqref{eq:tilde f} for the second inequality.
One now establishes that $t_1=T$ if $S$ is chosen large enough. Let $t\in[t_0,t_1)$ and let $i\in\cP$, $j\in\cQ.$ Equations~\eqref{peter pan} and~\eqref{bambi} together gives
\begin{equation}\label{piou piou}\begin{split}
\big|x_i(t)-x_j(t)\big|&\geq \big|x_i(t_0)-x_j(t_0)\big|-\Big|\big(x_i(t)-B_\cP(t)\big)-\big(x_i(t_0)-B_\cP(t_0)\big)\Big|-\Big|B_\cP(t)-B_\cP(t_0)\Big|\\&\qquad-\Big|\big(x_j(t)-B_\cQ(t)\big)-\big(x_j(t_0)-B_\cQ(t_0)\big)\Big|-\Big|B_\cQ(t)-B_\cQ(t_0)\Big|\\
&\geq\frac{S}{N-1}-2\frac{a}{A_0}C\Big(N-1,\,\max_k\|f_k\|_{L^\infty}+a\,\overline{G'},\,T,\,a,\,\frac{1}{A_0},\,\overline{G'}\Big)\\
&\qquad-2\frac{a}{A_0}T\Big(\max_k\|f_k\|_{L^\infty}+a\,\overline{G'}\Big).
\end{split}
\end{equation}
If one chooses $S$ big enough so that the last expression above is larger than $2$, then the definition of $t_1$ implies $t_1=T$. For that purpose we set the parameter $S$ equal to:
\begin{equation}\label{def:S_N}\begin{split}
S_N:=2(N-1)+2(N-1)\frac{a}{A_0}C\Big(N-1,\,\max_k\|f_k\|_{L^\infty}+a\,\overline{G'},\,T,\,a,\,\frac{1}{A_0},\,\overline{G'}\Big)\\+2(N-1)\frac{a}{A_0}T\Big(\max_k\|f_k\|_{L^\infty}+a\,\overline{G'}\Big).\end{split}
\end{equation}
Remark that the constraint $S_N\geq N$ previously required in the case $N\geq 2$ is indeed satisfied with such a definition.

It is now possible to establish the announced estimate. First, in the case $0\leq t\leq t_0$, then by definition of $t_0$,
\begin{equation}\label{chante}
\forall\;i,j=1\dots N,\qquad\big|\big(x_i(t)-x_j(t)\big)-\big(x_i(0)-x_j(0)\big)\big|\leq 2S_N.
\end{equation}
Otherwise if $t_0\leq t<t_1=T$,
\begin{equation}
\begin{split}
&\forall\;i,j=1\dots N,\qquad\Big|\big(x_i(t)-x_j(t)\big)-\big(x_i(0)-x_j(0)\big)\Big|\\&\qquad\leq\Big|\big(x_i(t_0)-x_j(t_0)\big)-\big(x_i(0)-x_j(0)\big)\Big|+\big|x_i(t)-x_i(t_0)\big|+\big|x_j(t)-x_j(t_0)\big|
\end{split}
\end{equation}
The first term above is estimated similarly as~\eqref{chante}. For the other terms, suppose for instance that $i\in\cP$, then one uses again~\eqref{peter pan} and~\eqref{bambi} to get
\begin{equation}\label{danse}\begin{split}
\Big|x_i(t)-x_i(t_0)\Big|&\leq\Big|\big(x_i(t)-B_\cP(t)\big)-\big(x_i(t_0)-B_\cP(t_0)\big)\Big|+\Big|B_\cP(t)-B_\cP(t_0)\Big|\\
&\leq\frac{a}{A_0}C\Big(N-1,\,\max_k\|f_k\|_{L^\infty}+a\,\overline{G'},\,T,\,a,\,\frac{1}{A_0},\,\overline{G'}\Big)\\
&\qquad+\frac{a}{A_0}T\Big(\max_k\|f_k\|_{L^\infty}+a\,\overline{G'}\Big)\\
&=\frac{S_N-(N+1)}{2(N-1)}\leq \frac{S_N}{2}
\end{split}
\end{equation}
Thus, gathering the two estimates~\eqref{chante} and~\eqref{danse},
\begin{equation}\label{bio}
\forall\;i,j,\quad\forall\;t\in[0,T),\qquad\big|\big(x_i(t)-x_j(t)\big)-\big(x_i(0)-x_j(0)\big)\big|\leq 3S_N.
\end{equation}
In view of the definition of $S_N$ at~\eqref{def:S_N}, it is possible to define the function $C(N,\dots)$ such that
\begin{equation}
C\Big(N,\,\max_k\|f_k\|_{L^\infty},\,T,\,a,\,\frac{1}{A_0},\,\overline{G'}\Big)=3S_N
\end{equation}
It is a direct computation to check that the function $C$ is increasing with respect to any of its variables and~\eqref{bio} corresponds exactly to~\eqref{eq:announced estimate}.
\end{proof}

\subsubsection{Proof of Corollary~\ref{lem:non collapse}}
Since the $a_i$ are all positive, the non-neutral clusters hypothesis~\eqref{eq:no null partial sum} is satisfied and then, as a consequence of Theorem~\ref{thrm:borne uniforme} the trajectories are bounded by a constant $C$. Thus,
\begin{equation}
\sum_{i\neq j}a_i\,a_j\,G\big(|x_i(t)-x_j(t)|\big)\geq a_{i_0}\,a_{j_0}\,G\big(|x_{i_0}(t)-x_{j_0}(t)|\big)+\Big(\sum_{\substack{i\neq j\\\{i,j\}\neq\{i_0,j_0\}}}a_i\,a_j\Big)\,\min_{0<r\leq C}G(r).
\end{equation}
where $i_0\neq j_0$ are two fixed indices. The left-hand side in the inequality above is a constant of the motion~\eqref{lem:preservation Hamiltonien}. Since by hypothesis $G(r)\to+\infty$ as $r\to0^+$, we conclude that $$\inf_{t\in[0,T)}|x_{i_0}(t)-x_{j_0}(t)|>0.$$\qed

\subsubsection{Proof of Proposition~\ref{thrm:uniform relative bound}}
Let $S>0$ parameter very large fixed later. Suppose toward
a contradiction that there exists an initial datum $X\in\RR^{2N}\setminus\mC_T$, two indices $i_0$ and $j_0$ and a time $t_0$ such that
\begin{equation}
|x_{i_0}(t_0)-x_{j_0}(t_0)|\;\geq\;\max\limits_{i,j}|x_i-x_j|+S.
\end{equation}
By continuity of the trajectories, there exists an interval of time $[t_1,t_2]$ such that
\begin{equation}\label{Barcelona}
|x_{i_0}(t_1)-x_{j_0}(t_1)|=\max\limits_{i,j}|x_i-x_j|+\frac{S}{2},\qquad|x_{i_0}(t_2)-x_{j_0}(t_2)|=\max\limits_{i,j}|x_i-x_j|+S,
\end{equation}
and such that
\begin{equation}
\forall\;t\in[t_1,t_2],\quad|x_{i_0}(t)-x_{j_0}(t)|\geq\max\limits_{i,j}|x_i-x_j|+\frac{S}{2}.
\end{equation}
As a consequence of this last property and similarly as in the proof of Theorem~\ref{thrm:borne uniforme}, we can split the system of vortices into two non-empty 
subsets $P$ and $Q$ with $P\cup Q=\{1\dots N\}$ such that for all $t\in[t_1,t_2]$,
\begin{equation}\label{Madrid}
\min\limits_{i\in P}\min\limits_{j\in Q}|x_i(t)-x_j(t)|\geq\frac{S}{2N}.
\end{equation}
It can be assumed for instance that $i_0\in P$ and $j_0\in Q$.
Therefore, the vortices in the set $P$ evolves according to equation
\begin{equation}\label{Sevilla}
\frac{d}{dt}x_i(t)=\sum_{\substack{j\in P\\j\neq i}}a_j\nabla^\perp G\big(|x_i^N(t)-x_j^N(t)|\big)+F_i(x_i^N(t),t).
\end{equation}
where the external fields $F_i$ is the interaction with the vortices 
that belongs to $Q$. As a consequence of~\eqref{Madrid}, and if $S$ is chosen large enough, this field satisfies 
\begin{equation}\label{Malaga}
|F_i(x,t)|\;\leq\;\sup_{r\geq 1}\Big|\frac{d}{dr}G(r)\Big|.
\end{equation}
The analogous equation holds for the vortices of the set $Q$. 
We now want to apply Proposition~\ref{prop:borne reformule} to the dynamic of the 
cluster $P$ on $[t_1,t_2]$ given by~\eqref{Sevilla}. The non-neutral 
clusters hypothesis is satisfied for the cluster $P$ as a consequence 
of~\eqref{eq:no null sub partial sum} because $P$ is a strict subset of $\{1\dots N\}$. We obtain from Proposition~\ref{prop:borne reformule} a constant $C$ such that the dynamic of the cluster $P$ without external field $F_i$ is bounded by $C$. If we add the smooth external field $F_i(x,t)$ and since the bound given by Proposition~\ref{prop:borne reformule} is uniform, we end up with
\begin{equation}
\sup\limits_{t\in[t_1,t_2]}\;\max\limits_{i\in P}\;|x_i(t)-x_i(t_1)|\leq C+\int_{t_1}^{t_2}\sup_{x\in\RR^{2N}}|F_i(x,t')|\,dt'.
\end{equation}
Combining this with~\eqref{Malaga} gives
\begin{equation}\label{Zaragossa}
\sup\limits_{t\in[t_1,t_2]}\;\max\limits_{i\in P}\;|x_i(t)-x_i(t_1)|\leq C+T\,\sup_{r\geq 1}\Big|\frac{d}{dr}G(r)\Big|.
\end{equation}
Doing an analogous argument on the cluster $Q$ gives
\begin{equation}\label{Murcia}
\sup\limits_{t\in[t_1,t_2]}\;\max\limits_{i\in Q}|x_{i}(t)-x_{i}(t_1)|\leq C+T\,\sup_{r\geq 1}\Big|\frac{d}{dr}G(r)\Big|.
\end{equation}
On the other hand, the condition~\eqref{Barcelona} implies that
\begin{equation}\label{Valencia}
|x_{i_0}(t_2)-x_{i_0}(t_1)|\geq\frac{S}{4}\qquad\text{or}\qquad|x_{j_0}(t_2)-x_{j_0}(t_1)|\geq\frac{S}{4}.
\end{equation}
The two bounds~\eqref{Zaragossa} and~\eqref{Murcia} are in contradiction with~\eqref{Valencia} if $S$ is chosen large enough.
\qed

\subsection{Proofs of the lemmas for Theorem~\ref{thrm:generalization of no partial sum with ordinary}}\label{sec:proofs 2}

\subsubsection{Proof of Lemma~\ref{lem:Dirac}}
Let $P_X(t)$ be the distribution defined by~\eqref{def:Dirac sum}, 
assumed well-defined for all $t\in[0,T[$. Let $\varphi\in\cD(\RR^2)$. Then,
\begin{equation}
\frac{d}{dt}\big<P_X(t),\varphi\big>_{\cD'\cD}
=\frac{d}{dt}\sum_{i=1}^Na_i\,\varphi\big(x_i(t)\big)=\sum_{i=1}^Na_i\nabla\varphi\big(x_i(t)\big)\cdot\frac{d}{dt}x_i(t).
\end{equation}
The equations of motion~\eqref{eq:vortex equation} give
\begin{equation}\begin{split}
&\frac{d}{dt}\big<P_X(t),\varphi\big>_{\cD'\cD}=\sum_{i\neq j}a_ia_j
\nabla\varphi\big(x_i(t)\big)\cdot\frac{\big(x_j(t)-x_i(t)\big)^\perp}{|x_j(t)-x_i(t)|^2}\\
&=\frac{1}{2}\sum_{i\neq j}a_ia_j\bigg(\nabla\varphi\big(x_i(t)\big)
-\nabla\varphi\big(x_j(t)\big)\bigg)\cdot\frac{\big(x_j(t)-x_i(t)\big)^\perp}{|x_j(t)-x_i(t)|^2}.
\end{split}\end{equation}
Therefore,
\begin{equation}
\bigg|\frac{d}{dt}\big<P_X(t),\varphi\big>_{\cD'\cD}\bigg|\leq
\frac{1}{2}\sum_{i\neq j}|a_i a_j|\,\big\|\nabla^2\varphi\big\|_\infty.
\end{equation}
Then, $t\mapsto\big<P_X(t),\varphi\big>_{\cD'\cD}$ is Lipschitz and 
converges as $t\to T^-$. 
Since this holds for any $\varphi\in\cD(\RR^2)$, this implies that 
$P_X(t)$ converges in the sense of distributions towards some $P_X\in\cD'(\RR^2)$ as $t\to T^-$. 
There remains to prove that $P_X$ is actually a measure that takes 
the form given by~\eqref{eq:Dirac limit}. 
Consider now an increasing sequence $(t_n)$ converging towards $T^-$. 
We remark first that it is always possible, up to an omitted extraction 
of the sequence, to reduce the problem to
\begin{equation}\label{alternative}
\text{either}\qquad x_i(t_n)\longrightarrow x_i^\ast\qquad\text{or}\qquad\big|x_i(t_n)\big|\longrightarrow+\infty
\end{equation} 
for some $X^\ast\in\RR^{2N}$. Indeed, if 
$\big|x_i(t_n)\big|\longrightarrow+\infty$ is not satisfied then 
there exists an extraction such that $x_i(t_n)$ stays bounded. But 
if it stays bounded then another extraction makes this sequence 
converge towards some $x_i^\ast$. Repeating this process for $i$ from 
$1$ to $N$ gives~\eqref{alternative}. Now that~\eqref{alternative} holds, define
\begin{equation}
b_i=\left\{\begin{array}{ll}0&\quad\text{ if }\big|x_i(t_n)\big|\longrightarrow+\infty,\\ 1&\quad\text{ either.}\end{array}\right.
\end{equation}
Therefore it holds
\begin{equation}
\sum_{i=1}^Na_i\,\delta_{x_i(t_n)}\longrightarrow\sum_{i=1}^Na_i\,b_i\,\delta_{x_i^\ast}\qquad\text{as}\;n\to+\infty,
\end{equation}
in the distributionnal sense.
By uniqueness of the limit, it is possible to identify
\begin{equation}
P_X = \sum_{i=1}^Na_i\,b_i\,\delta_{x_i^\ast}.
\end{equation}
The fact that the convergence of $P_X(t)$ towards $P_X$ in $\cD'$ is 
actually a convergence in the weak sense of measure comes from the 
fact that the measure $P_X(t)$ is bounded by $\sum_i|a_i|$ for all $t$.\qed

\subsubsection{Proof of Theorem~\ref{thrm:generalization of no partial sum with ordinary}}
\textbf{$\bullet\;$ Step 1: }Consider the $X^\ast\in\RR^{2N}$ given 
by Lemma~\ref{lem:Dirac}. Let $z\in\RR^2$ such that for all $i=1\dots N$, $z\neq x^\ast_i$. 
We are going to prove that for all $i=1\dots N$,
\begin{equation}\label{Ra}
\liminf\limits_{t\to T^-}|x_i(t)-z|>0.
\end{equation}
Suppose toward
a contradiction that there exists $A\subseteq\{1\dots N\}$ 
with $A\neq\emptyset$ such that for all $i\in A$,
\begin{equation}\label{Khepri}
\liminf\limits_{t\to T^-}|x_i(t)-z|=0.
\end{equation}
This set $A$ can be chosen such that for all $i\notin A$,
\begin{equation}\label{Nout}
\liminf\limits_{t\to T^-}|x_i(t)-z|>0.
\end{equation}
Define
\begin{equation}\label{Hapi}\begin{split}
&d_z^1:=\min\big\{|x_i^\ast-z|:i=1\dots N\big\}>0,\\
&d_z^2:=\min\big\{\liminf\limits_{t\to T^-}|x_i(t)-z|:i\notin A\big\}>0,\\
&d^\ast_z:=\min\big\{d_z^1,\;d_z^2\big\}>0,
\end{split}
\end{equation}
where by convention, for the definition of $d_z^2$, the minimum of 
the empty set is $+\infty$.
Let $\varphi$ be a $\cC^\infty$ function supported on the ball $\cB(z,d^\ast_z/2)$ and equal to $1$ on the ball $\cB(z,d^\ast_z/4)$. 
As a consequence of Lemma~\ref{lem:Dirac} (and Theorem~\ref{thrm:borne uniforme} to obtain that $b_i=1$ for all $i$) and by 
definition of $d^\ast_z$ we have,
\begin{equation}\label{Re}
\Big<\sum_{i=1}^Na_i\delta_{x_i(t)},\;\varphi\Big>_{\cD',\cD}\;\longrightarrow
\;\Big<\sum_{i=1}^Na_i\delta_{x_i^\ast},\;\varphi\Big>_{\cD',\cD}= 0\qquad\text{as }t\to T^-.\end{equation}
Using now~\eqref{Khepri} and~\eqref{Nout}, we infer the existence of 
an increasing sequence $(t_n)_{n\in\NN}$ converging towards $T^-$ 
such that for all $n\in\NN$,
\begin{equation}
\forall\;i\in A,\quad x_i(t_n)\in\cB\Big(z,\,\frac{d^\ast_z}{4}\Big)
\qquad\mathrm{and}\qquad\forall\;i\notin A,\quad x_i(t_n)\notin\cB\Big(z,\,\frac{d^\ast_z}{2}\Big).
\end{equation}
With the definition of $\varphi$ and the definition of $d_z^\ast$ given at~\eqref{Hapi}, holds for all $n\in\NN$,
\begin{equation}\label{Amon}
\Big<\sum_{i=1}^Na_i\delta_{x_i(t_n)},\;\varphi\Big>_{\cD',\cD}=\sum_{i\in A}a_i.
\end{equation}
As a consequence of the non-neutral clusters hypothesis~\eqref{eq:no null partial sum} 
we have $\sum_{i\in A}a_i\neq 0$. Therefore,  Equations~\eqref{Re} and~\eqref{Amon} are in contradiction and then~\eqref{Ra} holds.

\textbf{$\bullet\;$ Step 2: }
Assume now that a given vortex $x_i(t)$ has two adherence points. 
By~\eqref{Ra}, these two points must be some $x_i^\ast$. For 
instance, $x_j^\ast$ and $x_k^\ast$ with $j$ and $k$ such that $x_j^\ast\neq x_k^\ast$. Define the smallest distance between $x_j^\ast$ 
and any other possible adherence point
\begin{equation}
r^\ast_j:=\min\Big\{r>0:\exists\;l=1\dots N,\;|x_j^\ast-x_l^\ast|=r\Big\}.
\end{equation}
Consider then the circle
\begin{equation}
\cS:=\Big\{x\in\RR^2:|x-x_j^\ast(t)|=\frac{1}{2}r^\ast_j\Big\}.
\end{equation}
Since $x_j^\ast$ is inside the ball of radius $r^\ast_j/2$ and 
$x_k^\ast$ is outside, since these two points are adherence points 
for the dynamics of $x_i(t)$ as $t\to T^-$ and since the 
trajectories are continuous, there exist an increasing sequence of time $(t_n)_{n\in\NN}$ converging towards $T^-$ such that
\begin{equation}
\forall\;n\in\NN,\quad x_i(t_n)\in\cS.
\end{equation}
By compactness of $\cS$, it can be assumed that, up to an extraction, $x_i(t_n)\to x^\ast\in\cS$ as $n\to\infty$.
By definition of $r^\ast_j$, for all $l=1\dots N$, $x_l^\ast\notin\cS$. 
These two facts together are in contradiction with~\eqref{Ra} and 
this concludes the proof of Theorem~\ref{thrm:generalization of no partial sum with ordinary}.
\qed

\subsection{Proofs of the lemmas for Theorem~\ref{thrm:Improved Marchioro Pulvirenti}}\label{sec:proofs 3}

\subsubsection{Proof of Lemma~\ref{lem:reformulation}}
First, the conclusion of Theorem~\ref{thrm:Improved Marchioro Pulvirenti} can be 
formulated as follows
\begin{equation}\label{Sekhmet}
\cL^{2N}\Big\{X\in\RR^{2N}:\exists\;T_X\in\RR_+,
\quad\liminf\limits_{t\to T_X^-}\;\min\limits_{i\neq j}\big|x_i(t)-x_j(t)\big|=0\Big\}\;=0,
\end{equation}
where $\cL^d$ refers to the Lebesgue measure of dimension $d$.
It is possible to reduce the problem to bounded intervals of time 
and bounded regions of space by rewriting~\eqref{Sekhmet} as follows.
\begin{equation}\label{Sekhmet2}\begin{split}
\cL^{2N}\bigg(\bigcup_{T=1}^{+\infty}\bigcup_{\rho=1}^{+\infty}\Big\{X\in\RR^{2N}:\;&
\exists\;T_X\in[0,T],
\quad\liminf\limits_{t\to T_X^-}\;\min\limits_{i\neq j}\big|x_i(t)-x_j(t)\big|=0\\
&\text{and}\qquad
\;\max_{i\neq j}\big|x_i(t=0)-x_j(t=0)\big|\leq\rho
\Big\}\bigg)=0.
\end{split}
\end{equation}
Indeed, we can directly check that the two sets appearing respectively in~\eqref{Sekhmet} and~\eqref{Sekhmet2} are equal. 
Since the reunion in~\eqref{Sekhmet2} is a countable reunion, then to conclude that~\eqref{Sekhmet} holds it is enough to prove that for all $T>0$ and $\rho>0$,
\begin{equation}\label{Horus}\begin{split}
\cL^{2N}\Big\{X\in\RR^{2N}:\;&
\exists\;T_X\in[0,T],
\quad\liminf\limits_{t\to T_X^-}\;\min\limits_{i\neq j}\big|x_i(t)-x_j(t)\big|=0\\
&\text{and}\qquad
\;\max_{i\neq j}\big|x_i(t=0)-x_j(t=0)\big|\leq\rho
\Big\}=0.\end{split}
\end{equation}

Now, let $T>0$ and $\rho>0$. For $i$ fixed in $\{1\dots N\}$ denote 
by $\cT_i$ the isomorphism that gives the position of the point 
vortices $(x_k)_{k=1}^N$ knowing the position of $x_i$ and knowing 
the differences $(y_{ij})_{j\neq i}$. In other words define,
\begin{equation}
\begin{array}{cccc}
\cT_i:&\RR^2\times\RR^{2(N-1)}&\to&\RR^{2N}\\
\;&(x,Y)&\mapsto&(x+y_{i1}\dots x+y_{i(i-1)},\;x,\;x+y_{i(i+1)}\dots x+y_{iN}).
\end{array}
\end{equation}
Thus,
\begin{equation}\label{Isis}
\begin{split}
&\bigg\{X\in\RR^{2N}:\exists\;T_X\in[0,T],\quad
\liminf\limits_{t\to T_X^-}\;\min\limits_{\substack{j=1\dots N\\j\neq i}}\;\big|x_i(t)-x_j(t)\big|=0\quad\text{and}\quad\max_{i\neq j}\big|x_i(t=0)-x_j(t=0)\big|\leq\rho.\bigg\}\\&=\bigcap_{i=1}^N\cT_i\bigg[\RR^2\times\Big\{Y_i:=(y_{ij})_{j\neq i}\in\cB(0,\rho)^{2(N-1)}:
\exists\;T_X\in[0,T],\quad\liminf_{t\to T_X^-}\;\min\limits_{j\neq i}\big|y_{ij}(t)\big|=0\Big\}\bigg].
\end{split}
\end{equation}
Using now hypothesis~\eqref{eq:the set with measure zero} and the Fubini theorem,
\begin{equation}\label{Osiris}
\cL^{2N}\bigg(\RR^2\times\Big\{Y_i:=(y_{ij})_{j\neq i}\in\cB(0,\rho)^{2(N-1)}:
\exists\;T_X\in[0,T],\quad\liminf_{t\to T_X^-}\;\min\limits_{j\neq i}\big|y_{ij}(t)\big|=0\Big\}\bigg)=0.
\end{equation}
Since $\cT_i$ is a linear map, it is absolutely continuous and 
therefore maps any sets of Lebesgue measure $0$ into sets of 
Lebesgue measure $0$.
Therefore, combining this fact with~\eqref{Isis} and~\eqref{Osiris} gives~\eqref{Horus}.\qed

\subsubsection{Proof of Lemma~\ref{lem:reformulation convergence}}
Let $T>0$, $\rho>0$ and $\varepsilon>0$. For $i\neq j$, we define 
the set of initial datum such that occurs an $\varepsilon$-collapse 
between the two vortices $x_i$ and $x_j$
\begin{equation}\label{Seth}
\Gamma_{ij}^{\varepsilon,\rho}:=\Big\{Y_i=(y_{il})_{l\neq i}\in\cB(0,\rho)^{2(N-1)}:
\;\exists\;t\in[0,T],\quad\big|y_{ij}(t)\big|\leq\varepsilon\Big\}.
\end{equation}
We also define the time at which occurs the $\varepsilon$-collapse. Let $Y_i\in\bigcup_{j\neq i}\Gamma_{ij}^{\varepsilon,\rho}$,
\begin{equation}
T_{Y_i}^\varepsilon:=\inf\Big\{t\in[0,T]:\min\limits_{j\neq i}\;\big|y_{ij}(t)\big|\leq\varepsilon\Big\}.
\end{equation}
We are also interested in the situations where other collapses 
occur, far from $x_i$.  This corresponds to the $\varepsilon$-collapses of vector $y_{jk}:=x_j-x_k=y_{ij}-y_{ik}$. Let $k\neq i,j$, define
\begin{equation}\label{Anubis}
\Gamma_{ijk}^{\varepsilon,\rho}
:=\Big\{Y_i\in\Gamma_{ij}^{\varepsilon,\rho}:\exists\;t<T_{Y_i}^\varepsilon,\quad\big|y_{ij}(t)-y_{ik}(t)\big|\leq\varepsilon\Big\}.
\end{equation}
The fact that $\Gamma_{ijk}^{\varepsilon,\rho}$ gives information on 
whether another $\varepsilon$-collapse occurs far from $x_i$ with 
$x_j$ before the expected $\varepsilon$-collapse between $x_i$ and 
$x_j$ implies the following inclusion.
\begin{equation}\label{Maat}
\Gamma_{ijk}^{\varepsilon,\rho}\;\subseteq\;\Gamma_{kj}^{\varepsilon,2\rho}\setminus\Gamma_{kji}^{\varepsilon,2\rho}.
\end{equation}
This inclusion must be understood as follows. If occurs an $\varepsilon$-collapse between $x_j$ and $x_k$ before the first $\varepsilon$-collapse between $x_i$ and $x_j$ (left-hand side of the inclusion above), then in particular we have an $\varepsilon$-collapse between $x_j$ and $x_k$ (the set $\Gamma_{kj}^{\varepsilon,2\rho}$ in the right-hand side above). Yet, since we do not have an $\varepsilon$-collapse between $x_j$ and $x_i$ before the $\varepsilon$-collapse between $x_j$ and $x_k$, we can remove the set $\Gamma_{kji}^{\varepsilon,2\rho}$ in the right-hand side of the inclusion above.
Another important inclusion is
\begin{equation}\label{Nekhbet}
\Gamma_{ijk}^{\varepsilon,\rho}\;\subseteq\;\Gamma_{ijk}^{\varepsilon,2\rho}.
\end{equation}
These two definitions~\eqref{Seth} and~\eqref{Anubis} study the 
$\varepsilon$-collapse on the exact system with kernel $G_s$. We 
need the same definitions with the regularized kernels $G_{s,\varepsilon}$.
\begin{equation}\label{Thot}
\begin{split}
&\widehat{\Gamma}_{ij}^{\varepsilon,\rho}:=\Big\{Y_i=(y_{il})_{l\neq i}\in\cB(0,\rho)^{2(N-1)}:
\;\exists\;t\in[0,T],\quad\big|y_{ij}^\varepsilon(t)\big|\leq\varepsilon\Big\},\\
&\widehat{T}_{Y_i}^\varepsilon:=\inf\Big\{t\in[0,T]:\min\limits_{j\neq i}\;\big|y_{ij}^\varepsilon(t)\big|\leq\varepsilon\Big\},\\
&\widehat{\Gamma}_{ijk}^{\varepsilon,\rho}
:=\Big\{Y_i\in\widehat{\Gamma}_{ij}^{\varepsilon,\rho}:\exists\;t<\widehat{T}_{Y_i}^\varepsilon,\quad\big|y_{ij}^\varepsilon(t)-y_{ik}^\varepsilon(t)\big|\leq\varepsilon\Big\}.
\end{split}
\end{equation}
One remarks now that as long as the quantities $\big|y_{ij}\big|$ and 
$\big|y_{ij}-y_{ik}\big|$ remain higher than $\varepsilon$ for all 
$j\neq i$ and $k\neq i,j$, then the dynamics of $y_{ij}$ and $y_{ij}^\varepsilon$ 
coincide as a consequence of~\eqref{eq:condition on G epsilon 1}. 
This property implies in particular, using the sets defined at~\eqref{Seth},~\eqref{Anubis} and~\eqref{Thot}, 
\begin{equation}\label{Hator}
\Gamma_{ij}^{\varepsilon,\rho}\setminus\bigg(\bigcup_{k\neq i,j}\Gamma_{ijk}^{\varepsilon,\rho}\bigg)
=\widehat{\Gamma}_{ij}^{\varepsilon,\rho}\setminus\bigg(\bigcup_{k\neq i,j}\widehat{\Gamma}_{ijk}^{\varepsilon,\rho}\bigg).
\end{equation}
The hypothesis of Lemma~\ref{lem:reformulation convergence} can be rephrased as follows: for all $\rho>0$, and for all $i\neq j$,
\begin{equation}\label{Geb}
\cL^{2N}\Big(\widehat{\Gamma}_{ij}^{\varepsilon,\rho}\Big)\longrightarrow 0\qquad\text{as }\varepsilon\to0^+.
\end{equation}
Concerning the conclusion, it is enough to prove that for all $\rho>0$, and for all $i\neq j$,
\begin{equation}\label{Khnoum}
\cL^{2N}\Big(\Gamma_{ij}^{\varepsilon,\rho}\Big)\longrightarrow 0\qquad\text{as }\varepsilon\to0^+,
\end{equation}
because for all $i=1\dots N$ the following equality holds:
\begin{equation}\label{Shou}
\Big\{Y_i:=(y_{ij})_{j\neq i}\in\cB(0,\rho)^{2(N-1)}:\exists\;T_X\in[0,T],
\quad\liminf_{t\to T_X^-}\;\min\limits_{j\neq i}\big|y_{ij}(t)\big|=0\Big\}
=\bigcap_{n=1}^{+\infty}\bigcup_{j\neq i}\Gamma_{ij}^{\frac{1}{n},\rho}.
\end{equation}
The fact that the convergences~\eqref{Geb} imply the 
convergences~\eqref{Khnoum} is given by the following computations. First, using~\eqref{Maat} we get
\begin{equation}\label{Papyrus}
\Gamma_{ij}^{\varepsilon,\rho}=\Bigg[\Gamma_{ij}^{\varepsilon,\rho}
\setminus\bigg(\bigcup_{k\neq i,j}\Gamma_{ijk}^{\varepsilon,\rho}\bigg)\Bigg]
\cup\bigg(\bigcup_{k\neq i,j}\Gamma_{ijk}^{\varepsilon,\rho}\bigg)\subseteq\Bigg[\Gamma_{ij}^{\varepsilon,\rho}\setminus
\bigg(\bigcup_{k\neq i,j}\Gamma_{ijk}^{\varepsilon,\rho}\bigg)\Bigg]\cup\bigg(\bigcup_{k\neq i,j}\Gamma_{kj}^{\varepsilon,2\rho}\setminus\Gamma_{kji}^{\varepsilon,2\rho}\bigg).
\end{equation}
One remarks that it is possible to do the same computation with the remaining term on the very right using again~\eqref{Maat} and this gives
\begin{equation}\begin{split}
\Gamma_{kj}^{\varepsilon,2\rho}\setminus\Gamma_{kji}^{\varepsilon,2\rho}
&=\Bigg[\Gamma_{kj}^{\varepsilon,2\rho}\setminus\bigg(\bigcup_{l\neq k,j}\Gamma_{kjl}^{\varepsilon,2\rho}\bigg)\Bigg]\cup\bigg(\bigcup_{l\neq i,j,k}\Gamma_{kjl}^{\varepsilon,2\rho}\bigg)\\
&\subseteq\Bigg[\Gamma_{kj}^{\varepsilon,2\rho}\setminus\bigg(\bigcup_{l\neq k,j}\Gamma_{kjl}^{\varepsilon,2\rho}\bigg)\Bigg]
\cup\bigg(\bigcup_{l\neq i,j,k}\Gamma_{lj}^{\varepsilon,4\rho}\setminus\Gamma_{ljk}^{\varepsilon,4\rho}\bigg).
\end{split}
\end{equation}
Doing the same computation as above recursively until the residual 
term is empty and using~\eqref{Nekhbet} transforms~\eqref{Papyrus} into
\begin{equation}
\Gamma_{ij}^{\varepsilon,\rho}\subseteq\bigcup_{k\neq j}
\Bigg[\Gamma_{kj}^{\varepsilon,\,2^N\!\!\rho}\setminus\bigg(\bigcup_{l\neq j,k}\Gamma_{kjl}^{\varepsilon,\,2^N\!\!\rho}\bigg)\Bigg].
\end{equation}
Using now~\eqref{Hator}, we finally get
\begin{equation}
\Gamma_{ij}^{\varepsilon,\rho}\subseteq\bigcup_{k\neq j}
\Bigg[\widehat{\Gamma}_{kj}^{\varepsilon,\,2^N\!\!\rho}\setminus
\bigg(\bigcup_{l\neq j,k}\widehat{\Gamma}_{kjl}^{\varepsilon,\,2^N\!\!\rho}\bigg)\Bigg]
\subseteq\bigcup_{k\neq j}\widehat{\Gamma}_{kj}^{\varepsilon,\,2^N\!\!\rho}.\end{equation}
Thus the convergences~\eqref{Geb} imply the convergences~\eqref{Khnoum} and the lemma is proved.\qed

\subsubsection{Proof of Lemma~\ref{lem:collapses}}
Let $i\in\{1\dots N\}$, $\varepsilon>0$ and $\rho>0$. 

\textbf{Step 1.}
Let $a>0$. We define a kernel $L_a$ by
\begin{equation}\label{Poseidon}
L_a(q):=q^{-2-a}.
\end{equation}
and we associate to the kernel $L_a$ its $\varepsilon$-regularization $L_{a,\varepsilon}$ as defined by~\eqref{eq:condition 
on G epsilon 1}-\eqref{eq:condition on G epsilon 4}. From this we 
define the function
\begin{equation}
\Phi(Y_i):=\sum_{j\neq i}L_{a,\varepsilon}\big(|y_{ij}|\big).
\end{equation}
This function is all the most high valued as the system is close to 
collapse with vortex $x_i$. Denote by $\mS^t_{i,\varepsilon}$ the flow of the modified 
system~\eqref{eq:evolution Y} with the regularized kernel $G_{s,\varepsilon}$. This gives
\begin{equation}\label{Hera}\begin{split}
&\frac{d}{dt}\Phi\Big(\mS^t_{i,\varepsilon} Y_i\Big)=\sum_{j\neq i}\nabla L_{a,\varepsilon}\big(|y_{ij}^\varepsilon(t)|\big)
\cdot\frac{d}{dt}y_{ij}^\varepsilon(t),\\
&=\sum_{j\neq i}\nabla L_{a,\varepsilon}\big(|y_{ij}^\varepsilon(t)|\big)
\cdot\!\bigg[(a_i+a_j)\nabla^\perp G_{s,\varepsilon}\big(|y_{ij}^\varepsilon|\big)+\!\sum_{k\neq i,j}\!a_k\Big(\nabla^\perp G_{s,\varepsilon}\big(|y_{ik}^\varepsilon|\big)
+\nabla^\perp G_{s,\varepsilon}\big(|y_{ij}^\varepsilon-y_{ik}^\varepsilon|\big)\Big)\bigg],\\
&=\sum_{j\neq i}\sum_{k\neq i,j}a_k\nabla L_{a,\varepsilon}\big(|y_{ij}^\varepsilon(t)|\big)
\cdot\Big(\nabla^\perp G_{s,\varepsilon}\big(|y_{ik}^\varepsilon|\big)
+\nabla^\perp G_{s,\varepsilon}\big(|y_{ij}^\varepsilon-y_{ik}^\varepsilon|\big)\Big),
\end{split}
\end{equation}
where for the last equality we used the identity $\nabla f\cdot\nabla^\perp g=0$ 
that holds for $f$ and $g$ two radial functions. Thus,
\begin{equation}\label{Demeter}
\bigg|\frac{d}{dt}\Phi\Big(\mS^t_{i,\varepsilon} Y_i\Big)\bigg|\leq\Psi\Big(\mS^t_{i,\varepsilon} Y_i\Big),
\end{equation}
where
\begin{equation}\label{Aphrodite}
\Psi\big(Y_i\big):=\sum_{j\neq i}\sum_{k\neq i,j}a_k
\big|\nabla L_{a,\varepsilon}\big(|y_{ij}|\big)\big|\;\Big(\big|\nabla G_{s,\varepsilon}
\big(|y_{ik}|\big)\big|+\big|\nabla G_{s,\varepsilon}\big(|y_{ij}-y_{ik}|\big)\big|\Big).
\end{equation}
We now observe, recalling the definition of $G_s$ given 
at~\eqref{def:Green functions}, that the $\varepsilon$-
regularization~\eqref{eq:condition on G epsilon 1}-\eqref{eq:condition on G epsilon 4} 
implies, by a direct computation using polar coordinates,
\begin{equation}\label{Artemis}
\int_{\cB(0,\rho)}\big|\nabla G_{s,\varepsilon}\big(|y|\big)|dy\leq C\left\{\begin{array}{ll}
1&\quad\text{if }s>0.5,\\
\log(1/\varepsilon)&\quad\text{if } s=0.5,\\
\varepsilon^{2s-1}&\quad\text{if } s<0.5,
\end{array}\right.
\end{equation}
where $\cB(0,\rho)$ is the Euclidean ball on $\RR^2$. T
he constant $C$ depends on $\rho$ and $s$. Similarly with the definition of the kernel $L_a$ at~\eqref{Poseidon} and since $a>0$,
\begin{equation}\label{Apollon}
\int_{\cB(0,\rho)}L_{a,\varepsilon}\big(|y|\big)dy
\leq C\varepsilon^{-a}\qquad\mathrm{and}
\qquad\int_{\cB(0,\rho)}\big|\nabla L_{a,\varepsilon}\big(|y|\big)|dy\leq C\varepsilon^{-1-a}.
\end{equation}
Therefore using~\eqref{Apollon} with the definition of $\Phi$ gives
\begin{equation}\label{Hermes}
\int_{\cB(0,\rho)^{N-1}}\Phi(Y_i)\,dY_i\leq C\varepsilon^{-a}.
\end{equation}
Similarly, the definition of $\Psi$ given at~\eqref{Aphrodite} gives
\begin{equation}\begin{split}
&\int_{\cB(0,\rho)^{N-1}}\Psi(Y_i)\,dY_i\\
&=\int_{\cB(0,\rho)^{N-1}}\Bigg[\sum_{j\neq i}\sum_{k\neq i,j}a_k
\big|\nabla L_{a,\varepsilon}\big(|y_{ij}|\big)\big|
\;\Big(\big|\nabla G_{s,\varepsilon}\big(|y_{ik}|\big)\big|
+\big|\nabla G_{s,\varepsilon}\big(|y_{ij}-y_{ik}|\big)\big|\Big)\Bigg]\prod_{\substack{l=1\\l\neq i}}^Ndy_{il}\\
&=2\Big(\sum_{j\neq i}\sum_{k\neq i,j}a_k\Big)\bigg(\int_{\cB(0,\rho)}dy\bigg)^{N-3}
\bigg(\int_{\cB(0,\rho)}\big|\nabla L_{a,\varepsilon}\big(|y|\big)\big|dy\bigg)
\bigg(\int_{\cB(0,\rho)}\big|\nabla G_{s,\varepsilon}\big(|y|\big)\big|dy\bigg),
\end{split}
\end{equation}
and then, using~\eqref{Artemis} and~\eqref{Apollon},
\begin{equation}\label{Dionysos}
\int_{\cB(0,\rho)^{N-1}}\Psi(Y_i)\,dY_i\leq C\varepsilon^{-2-a}\left\{\begin{array}{ll}
\varepsilon&\quad\text{if }s>0.5,\\
\varepsilon\log(1/\varepsilon)&\quad\text{if } s=0.5,\\
\varepsilon^{2s}&\quad\text{if } s<0.5.
\end{array}\right.
\end{equation}

\textbf{Step 2.} It is now possible to integrate $\Phi$ along the 
flow. We obtain
\begin{equation}
\begin{split}
\int_{\cB(0,\rho)^{N-1}}\sup_{t\in[0,T]}\Phi\big(\mS^t_{i,\varepsilon}Y_i\big)\,dY_i
\leq\int_{\cB(0,\rho)^{N-1}}\Phi(Y_i)\,dY_i
+\int_{\cB(0,\rho)^{N-1}}\int_0^T\bigg|\frac{d}{dt}\Phi\Big(\mS^t_{i,\varepsilon} Y_i\Big)\bigg|\,dt\,dY_i
\end{split}
\end{equation}
Using~\eqref{Demeter} in the estimate above gives
\begin{equation}\label{Hephaistos}
\int_{\cB(0,\rho)^{N-1}}\sup_{t\in[0,T]}\Phi\big(\mS^t_{i,\varepsilon}Y_i\big)\,dY_i
\leq\int_{\cB(0,\rho)^{N-1}}\Phi(Y_i)\,dY_i
+\int_{\cB(0,\rho)^{N-1}}\int_0^T\Psi\Big(\mS^t_{i,\varepsilon} Y_i\Big)\,dt\,dY_i.
\end{equation}
Using the Fubini theorem in~\eqref{Hephaistos} and the Liouville theorem~\ref{lem:Liouville Theorem for the modified dynamics} leads to
\begin{equation}\label{Circee}\begin{split}
\int_{\cB(0,\rho)^{N-1}}\sup_{t\in[0,T]}\Phi\big(\mS^t_{i,\varepsilon}Y_i\big)\,dY_i
&\leq\int_{\cB(0,\rho)^{N-1}}\Phi(Y_i)\,dY_i
+\int_0^T\int_{\mS^t_{i,\varepsilon}\cB(0,\rho)^{N-1}}\Psi\big( Y_i\big)\,d\mS^{-t}_{i,\varepsilon}Y_i\,dt\\
&=\int_{\cB(0,\rho)^{N-1}}\Phi(Y_i)\,dY_i
+\int_0^T\int_{\mS^t_{i,\varepsilon}\cB(0,\rho)^{N-1}}\Psi\big( Y_i\big)\,dY_i\,dt.
\end{split}
\end{equation}
We now make use of hypothesis~\eqref{eq:no null sub partial sum} on the intensities of the vorticies. Indeed, this hypothesis allows us to use Theorem~\ref{thrm:uniform relative bound} which states the existence of a constant $C'$ independent on $\varepsilon$ (but dependent on $\rho$, $T$, $s$ and the $a_i$) such that
\begin{equation}
\mS^t_{i,\varepsilon}\cB(0,\rho)^{N-1}\subseteq\cB(0,C')^{N-1}.
\end{equation}
Thus, the estimate~\eqref{Circee} above becomes
\begin{equation}\label{Hades}\begin{split}
\int_{\cB(0,\rho)^{N-1}}\sup_{t\in[0,T]}\Phi\big(\mS^t_{i,\varepsilon}Y_i\big)\,dY_i
&\leq\int_{\cB(0,\rho)^{N-1}}\Phi(Y_i)\,dY_i+\int_0^T\int_{\cB(0,C')^{N-1}}\Psi\big( Y_i\big)\,dY_i\,dt.\\
&\leq C\varepsilon^{-a}+C\,T\varepsilon^{-2-a}\left\{\begin{array}{ll}
\varepsilon&\quad\text{if }s>0.5,\\
\varepsilon\log(1/\varepsilon)&\quad\text{if } s=0.5,\\
\varepsilon^{2s}&\quad\text{if } s<0.5.
\end{array}\right.,
\end{split}
\end{equation}
where for the last estimate we used~\eqref{Hermes} and~\eqref{Dionysos}. 

\textbf{Step 3.} By definition of the function $\Phi$, there exists a constant $c>0$ such that,
\begin{equation}\begin{split}
&\Big\{Y_i=(y_{ij})_{j\neq i}\in \cB(0,\rho)^{N-1}:\min\limits_{j\neq i}\inf_{t\in[0,T]}\big|y_{ij}^\varepsilon(t)\big|
\leq\varepsilon\Big\}\\&\qquad\qquad\subseteq\Big\{Y_i=(y_{ij})_{j\neq i}\in\cB(0,\rho)^{N-1}
:\sup_{t\in[0,T]}\Phi\big(\mS^t_{i,\varepsilon}Y_i\big)\geq c\,\varepsilon^{-2-a}\Big\}.
\end{split}
\end{equation}
Combining this inclusion with the Bienaymé-Tchebycheff inequality gives
\begin{equation}\label{Hestia}\begin{split}
\cL^{2(N-1)}\Big\{Y_i=(y_{ij})_{j\neq i}\in \cB(0,\rho)^{N-1}:&
\min\limits_{j\neq i}\inf_{t\in[0,T]}\big|y_{ij}^\varepsilon(t)\big|
\leq\varepsilon\Big\}\\&\leq\frac{\varepsilon^{2+a}}{c}
\int_{\cB(0,\rho)^{N-1}}\sup_{t\in[0,T]}\Phi\big(\mS^t_{i,\varepsilon}Y_i\big)\,dY_i.\end{split}
\end{equation}
Using now~\eqref{Hades} in~\eqref{Hestia},
\begin{equation}
\cL^{2(N-1)}\Big\{Y_i=(y_{ij})_{j\neq i}\in \cB(0,\rho)^{N-1}:\min\limits_{j\neq i}\inf_{t\in[0,T]}\big|y_{ij}^\varepsilon(t)\big|
\leq\varepsilon\Big\}\leq\,C\left\{\begin{array}{ll}
\varepsilon&\quad\text{if }s>0.5,\\
\varepsilon\log(1/\varepsilon)&\quad\text{if } s=0.5,\\
\varepsilon^{2s}&\quad\text{if } s<0.5.
\end{array}\right.,
\end{equation}
where $C$ is a constant that depends on $T$, $\rho$, $N$, $s$ and on the $a_i$. The lemma is proved. \qed

\newpage
\vspace{1cm}
\noindent \textbf{\Large Acknowledgments}\vspace{0.3cm}

I would like to acknowledge my PhD advisors Philippe GRAVEJAT and 
Didier SMETS for their confidence and their scientific support, 
constructive criticism and suggestions during all my work on the point-vortex model. I would like also to acknowledge them for their meticulous rereadings and their advises during the redaction of this article.

The author acknowledges grants from the \emph{Agence Nationale de la Recherche} (ANR), for the project ``Ondes Dispersives 
Aléatoires" (ANR-18-CE40-0020-01). The problem considered in this 
paper was inspired from the Workshop of this ANR ODA on 27-28 June 
2019 in Laboratoire Paul Painlevé (Lille, France).

\bibliographystyle{plain}
\bibliography{bibliography}
\end{document}